\magnification=1200
\def\v{\varphi}
\def\tm{{\tilde M}}
\def\tn{{\tilde N}}

\def\diam{{\rm diam}}
\def\exp{{\rm exp}}

\def\f{{\cal F}}
\def\c{{\cal C}}
\def\g{{\cal G}}
\def\ch{{\cal H}}
\def\m{{\cal M}}
\def\p{{\cal P}}
\def\s{{\cal S}}
\def\j{{\cal J}}
\def\pl{{\parallel}}
\def\tg{{\tilde g}}

\def\tc{{\tilde c}}

\def\tp{{\tilde \Psi}}
\def\ep{\epsilon}
\def\op{\oplus}

\centerline{\bf Cocycles, symplectic structures and intersection}
\bigskip
\centerline{Ursula Hamenst\"adt}
\bigskip
\bigskip
\bigskip
\bigskip
\centerline{\bf 1. Introduction}
\bigskip
\bigskip
In the paper [B-C-G] Besson, Courtois and Gallot
proved the following 
remarkable theorem: Let $S$ be a closed rank 1 locally symmetric
space of noncompact type and let $M$ be a closed manifold of
negative curvature which is homotopy equivalent to $S$. If $S$ and
$M$ have the same volume and the same volume entropies (i.e. 
the same asymptotic growth rate of volumes of balls in their
universal covers), then $S$ and $M$ are isometric.

One application of this result is the solution of the so-called
{\sl conjugacy problem} for locally symmetric manifolds. This
problem can be stated in the following way: The {\sl marked length
spectrum} of a closed negatively curved manifold $(M,g)$ is the
function $\rho_M$ which assigns to each conjugacy class in the
fundamental group $\pi_1(M)$ of $M$ the length of the unique
$g$-geodesic representing the class. Two homotopy
equivalent
negatively curved  manifolds $M$, $N$ have {\sl the same marked length
spectrum} if there is an isomorphism $\Psi:\pi_1(M)\to \pi_1(N)$ such
that $\rho_N \circ \Psi= \rho_M$. This is equivalent to the
existence of a continuous {\sl time preserving conjugacy}
for the {\sl geodesic flows} for $M$ and $N$ ([H1]). 
Such a conjugacy is defined to be a homeomorphism 
$\Lambda$ of the unit tangent bundle 
$T^1M$ of $M$ onto the unit tangent bundle $T^1N$ of $N$
which is equivariant under the action of the geodesic flows
$\Phi^t$ on $T^1M$ and $T^1N$ ([H1]). The conjugacy problem now
asks whether a composition of $\Lambda$ with the time-t-map of the geodesic
flow on $T^1N$ for a suitable $t\in {\bf R}$ in the restriction
to $T^1M$ of the differential of an isometry of $M$ onto $N$.
This is known to be true for surfaces ([O]).

Since the volume entropy of a closed negatively curved manifold
equals the topological entropy of the geodesic flow, two manifolds
with time preserving conjugate geodesic flows have the same
volume entropy.

Let $X^0$ be the {\sl geodesic spray}, i.e. the generator of
the geodesic flow $\Phi^t$.
Recall
that there is a H\"older continuous $\Phi^t$-invariant decomposition
$TT^1 M = {\bf R} X^0 \op TW^{ss} \op TW^{su}$,
the so called {\sl Anosov splitting},
where $TW^{ss}$ (or $TW^{su}$) is
the tangent bundle of the {\sl strong stable (or strong unstable) foliation}
$W^{ss}$
(or $W^{su}$) on $T^1 M$.
The bundle
$TW^{ss}
\op TW^{su}$ is smooth, and
hence we obtain a smooth 1-form $\omega$ on $T^1 M$ by defining $\omega (X^0)
\equiv 1$ and $\omega (TW^{ss} \op TW^{su}) = 0$. This form $\omega$ is called
the {\sl canonical contact form}.
The differential form $\omega \wedge (d\omega)^{n-1}$ is the volume
form of the so called {\sl Sasaki metric} on $T^1M$, and the total
mass of $T^1M$ with respect to this volume form equals the
product of the volume of $M$ with the volume of the
$(n-1)$-dimensional unit sphere in ${\bf R}^n$.

The pull-back of the canonical contact form on $T^1N$ under a
time preserving conjugacy $\Lambda:T^1M\to T^1N$ of class $C^1$
is the canonical 
contact form on $T^1M$. This implies that $M$ and $N$ have the
same volume if their geodesic flows are 
$C^1$-time-preserving conjugate.
Therefore
the result of Besson, Courtois and Gallot
gives a positive answer to the conjugacy
problem if one of the manifolds under consideration is locally 
symmetric and if the conjugacy is assumed to be of class $C^1$.

One of the objectives of this note is to remove this regularity
assumption on the conjugacy. In Section
3 we show:
\bigskip
{\bf Theorem A:} {\it If $M$ and $N$are closed negatively curved
manifolds with the same marked length spectrum and if the Anosov
splitting of $TT^1M$ is of class $C^1$, then $M$ and $N$ have
the same volume.}
\bigskip
Manifolds with strictly 1/4- pinched sectional curvature
or locally symmetric manifolds 
have $C^1$ Anosov splitting. Thus
as a corollary from [B-C-G] we obtain:
\bigskip
{\bf Corollary:} {\it If $M$ has the same marked length spectrum
as a closed negatively curved locally symmetric space $S$, then
$M$ and $S$ are isometric.}
\bigskip
Recall that the {\sl space of geodesics} ${\cal G}\tm$ of the universal
cover $\tm$ of $(M,g)$ is the quotient of the unit tangent bundle
$T^1\tm$ of $\tm$ under the action of the geodesic flow $\Phi^t$.
Also ${\cal G}\tm$ can naturally be identified with
$\partial\tm\times \partial\tm -\Delta$ where $\partial \tm$ is
the {\sl ideal boundary} of $\tm$ and $\Delta$ is the diagonal.
The space of geodesics is a smooth manifold, and the differential
$d\omega$ of the canonical contact form projects to a smooth 
symplectic structure on ${\cal G}\tm$.
The fundamental group
$\pi_1 (M) $ of $M$ acts on ${\cal G}\tm$ as a group of 
symplectomorphisms.

The {\sl cross ratio} of the metric $g$ is a function on 
${\cal G}\tm \times {\cal G}\tm$. We show in Section 3 that this
function can be viewed as a coarse (or rather integrated) version of
the symplectic form $d\omega$ on ${\cal G}\tm$. 
This enables us to construct the measure
$\lambda_g$ 
on ${\cal G}\tm$ defined by the smooth volume form
$(d\omega)^{n-1}$ in a purely combinatorial way from the
cross ratio. The measure $\lambda_g$ is usually called the
{\sl Lebesgue Liouville current} of the metric $g$.
Our construction
is general, but unfortunately so far we are only able to show
that the measure
we obtain from it is exactly the Lebesgue Liouville
current (rather than some multiple of it) under the additional 
assumption that the Anosov splitting of $TT^1M$ is of class $C^1$.

A {\sl geodesic current}
for $M$ is a
locally finite Borel measure on 
the space of geodesics 
which is invariant under the action of the fundamental group and under
the exchange of the factors in the product decomposition of
${\cal G}\tm =\partial \tm\times \partial \tm-\Delta$.
If $M$ and $N$ are two homotopy equivalent closed negatively
curved manifolds, then there is a natural $\pi_1(M)=\pi_1(N)$-
equivariant H\"older homeomorphism of $\partial \tm$ onto $\partial \tn$
which induces an equivariant homoeomorphism
of ${\cal G}\tm $ onto ${\cal G}\tn$, in particular 
the spaces of geodesic currents for $M$ and $N$
are naturally identified. 

With this identification we show in Section 4:
\bigskip
{\bf Theorem B:} {\it Let $(M,g)$,
$(N,h)$ be homotopy equivalent closed
negatively curved manifolds.
If $\lambda_g=\lambda_h$ and if the Anosov splittings of
$TT^1M$ and $TT^1N$ are of class $C^1$, then
$M$ and $N$ have the same marked length spectrum.}
\bigskip

Recall
that
a H\"older
continuous {\sl additive cocycle} for the
geodesic flow $\Phi^t$ is a H\"older continuous function $\zeta \colon T^1 M
\times {\bf R} \to {\bf R}$ which satisfies $\zeta (v, s + t) = \zeta (v, s) + \zeta
(\Phi^s v, t)$ for all $v \in T^1 M$ and all $s, t \in {\bf R}$.
Two cocycles $\zeta, \xi$ are called {\sl cohomologous} if there is a H\"older
continuous function $\beta$ on $T^1 M$ such that $\zeta (v, t) - \xi (v,t)=
\beta (\Phi^t
v)  - \beta (v)$ for all $v \in T^1 M$ and all $t \in {\bf R}$. We
denote by $[\zeta]$ the cohomology class of $\zeta$.
Let $\f \colon
T^1 M \to T^1 M $ be the {\sl flip} $v \to \f v = - v;$ the flip $\f \zeta$
of a cocycle $\zeta$ is the cocycle $\f \zeta (v, t) = \zeta (\f \Phi^t v,
t)$. Two cocycles $\zeta, \xi$ are cohomologous if and only if this is true
for $\f \zeta$ and $\f\xi$. In other words, $\f$ induces an action on
cohomology classes of H\"older cocycles which we denote again by $\f$. We call
$\zeta$
{\sl quasi-invariant
under the flip} if the cocycles $\zeta$ and $\f \zeta$ are cohomologous, i.e.
if $[\zeta]$ is a fix point for the action of $\f$.

Consider now a closed negatively curved surface $M$.
Every flip invariant H\"older 
cohomology class for the geodesic flow defines
a unique cross ratio ([H5]), and in the two-dimensional case
the defining properties of a cross ratio ([H5]) show that 
we can obtain from a cross ratio a 
finitely additive signed measure on
${\cal G}\tm$ in a in a natural way.
In other words, there is a natural linear map which associates to
a H\"older cohomology class $[\zeta]$ a signed measure $\nu_\zeta$.
Now the {\sl length cocycle} $\ell$ of
the negatively curved manifold $(M,g)$ is defined
by $\ell(v,t)=t$ for all $v\in T^1M$ and all $t\in {\bf R}$. 
For the length cocycle $\ell$ of $g$ on the surface $M$,
$\nu_\ell$ is just
1/2 times the Lebesgue Liouville current of $g$.

On the other hand, if the H\"older
cohomology class is {\sl positive}, i.e. if it
can be represented by a cocycle $\bar \zeta$ which satisfies
$\bar \zeta (v,t)>0$ for all $v\in T^1M$ and all $t>0$, then
there is another more classical way to define a projective class
of currents on ${\cal G}\tm$. Namely, a 
suitable positive multiple $[a\zeta]$ of $[\zeta]$ (here $a>0$ is a 
constant depending on $[\zeta]$)
defines a {\sl Gibbs equilibrium state} which is a $\pi_1 (M)$-invariant
measure on  $\g \tm $, determined uniquely
up to a constant. This measure is a geodesic current
if and only if $[\zeta] = [\f\zeta]$, and if this is
satisfied we call it a
{\sl Gibbs-current of} $\zeta$ 
and denote its projective class in the space of
projectivized currents by $[\mu_\zeta]$.
The assignment $[\zeta] \to [\mu_\zeta]$ then is a map from a subset of the
projectivization of the space of flip invariant cohomology classes into the
projectivization of the space of geodesic currents.

If
$\ell$ is the length cocycle of the
Riemannian metric $g$ on $M$, then its Gibbs
current is the
{\sl Bowen-Margulis current} of $g$ which corresponds to the measure of
maximal entropy for the geodesic flow. 

Katok showed in [K] that the Lebesgue Liouville current and
the Bowen Margulis current of a negatively curved metric on a surface are
equivalent if and only if the metric has constant curvature.
With our above notation this result can be stated as follows:
If
$\ell$ is the length cocycle of a metric $g$ of negative
curvature on
the surface $M$ and if $\nu_\ell$ is equal to $\mu_\ell$ up
to a constant, then $g$ has constant curvature.

In Section 2 we obtain a generalization of Katok's result to 
arbitrary positive H\"older cohomology classes.
For its formulation denote by $[\nu]$ the class of a current
in the projectivization of the space of geodesic currents. We
show:
\bigskip

{\bf Theorem C:} {\it Let $M$ be a closed surface of genus $g \ge 2$. If
$[\mu_\ell] = [\nu_\ell]$ for a positive flip invariant H\"older
cohomology class
$\ell$,
then $\ell$ is the class of the length cocycle of a metric on $M$
of constant curvature. 
}
\bigskip

Before
we proceed we fix a few more notations. All our manifolds will be
closed and equipped with a fixed negatively curved metric. We denote by $T^1
M$ (or $ T^1 \tm$) the unit tangent bundle of $M$ (or $\tm$). There is a
natural $\pi_1 (M)$-equivariant projection $\pi \colon T^1 \tm \to \partial
\tm$. The pre-image of a point $\xi \in \partial \tm$ under $\pi$ equals the
stable manifold in $T^1 \tm$ of all directions pointing towards $\xi$. We let
$W^i (v)$ be the leaf of the foliation $W^i$ containing $v$
$(i={\rm ss},\ {\rm su})$ and write $P$ for
the canonical projections $T^1 M \to M$ and $T^1 \tm \to \tm$.

We will need the following facts which were pointed out
by Ledrappier ([L]):
The fundamental group $\pi_1 (M)$ of $M$ acts
naturally on 
the ideal boundary $\partial \tm$ 
of $\tm$ as a group of homeomorphisms. There
is a natural bijection between the space of H\"older cocycles for the
geodesic flow on $T^1M$ and the space of
H\"older cocycles for the action of $\pi_1(M)$ on 
$\partial \tm$. Such a H\"older 
cocycle for the action of $\pi_1(M)$ is a H\"older continuous
functions $\zeta_0 \colon \pi_1 (M) \times \partial \tm \to {\bf R}$ which
satisfies $\zeta_0 (\v_1 \v_2, \xi) = \zeta_0 (\v_1, \v_2 \xi) +
\zeta_0 (\v_2, \xi)$ for $\v_1, \v_2\in \pi_1 (M)$ and $\xi\in \partial \tm$.
Two such functions
$\zeta_0$, $\eta_0$ are {\sl cohomologous} if
there is a H\"older
continuous function $\beta_0$ on $\partial\tm$ such that $\zeta_0 (\v, \xi) -
\eta_0 (\v, \xi) = \beta_0 (\v \xi) - \beta_0 (\xi)$ for all $\v \in \pi_1
(M)$ and all $\xi \in \partial \tm$. If $\xi (\v)$
is the attracting fixed point for the action of $\v \in \pi_1 (M)$ on $\partial
\tm$, then $\zeta_0 (\v, \xi (\v))$ is
called the {\sl period} of $\zeta_0$ at $\v$.
Two H\"older cocycles are cohomologous if and only if they have the same
periods.

The natural map from H\"older cocycles for the geodesic
flow $\Phi^t$
to H\"older cocycles for the action of $\pi_1(M)$ on $\partial \tm$
preserves the equivalence relation which defines cohomology.
Moreover a H\"older class 
$[\zeta]$ for $\Phi^t$ is invariant under the flip if and only if
the corresponding class $[\zeta_0]$ for the action of $\pi_1(M)$
satisfies
$\zeta_0 (\v, \xi (\v)) = \zeta_0 (\v^{-1}, \xi(\v^{-1}))$ for all $\v \in
\pi_1 (M)$ (all this is discussed in [L]).
In particular, the space of H\"older classes for the geodesic flow of a 
negatively curved metric $g$ on $M$ as well as the notion of
flip-invariance and positivity do not depend on the choice of the metric.

Let $\m$ be the space of $\Phi^t$-invariant Borel probability measures on
$T^1 M$. For every H\"older cocycle $\zeta$ 
for $\Phi^t$and every
$\mu \in \m$ the integral $\int \zeta d\mu$ can be defined. This integral only
depends on the cohomology class $[\zeta]$ of $\zeta$.
More precisely, let
 $\alpha \colon T^1 M \to {\bf R}$ be a H\"older cocycle and write $f (v) =
\alpha (v, 1)$ for all $v \in T^1 M.$ Then $f$ is a H\"older
 continuous function
on $T^1 M$ and defines a H\"older cocycle $\alpha_f$ by the formula $\alpha_f
(v, t) = \int_0^t f (\Phi^s v)\; ds.$ If $v \in T^1 M$ is a periodic point for
$\Phi^t$ of period $\tau > 0$ then
$$\eqalign{\alpha_f (v, \tau) = & \int_0^\tau f (\Phi^s
v) \; ds = \int_0^\tau \alpha (\Phi^s v, 1) \; ds\cr
= \int^\tau_0 (\alpha (v, s + 1) - \alpha (v, s)) \; ds
=&\int^1_0 (\alpha (v, \tau + s) - \alpha (v, s)) \; ds
= \alpha (v, \tau)}$$
by the
cocycle equality for $\alpha$. But this just means that $\alpha$ and
$\alpha_f$ are equivalent\hfil\break
(compare [L]). Thus every H\"older cohomology class
can be represented by a H\"older continuous function $f$, and two functions $f$
and $g$ define the same H\"older cohomology class if and only if for every
element $\eta$ from the space $\m$ of $\Phi^t$- invariant Borel probability
measures on $T^1 M$ we have $\int f d\eta = \int g d\eta$ (this is the
Livshicv theorem for H\"older continuous functions).

Denote by $[f]$ the cohomology class of the cocycle $\alpha_f$ defined by $f$.
Recall that the {\sl flip} $\f \colon T^1 M \to T^1 M $ defined by $\f (v) =
-v$ is a diffeomorphism of $T^1 M$ which satisfies $\f \circ \Phi^t =
\Phi^{-t} \circ \f$. Call $[f]$ {\sl flip-invariant} if $[f]$ can be
represented by a function $g$ on $T^1 M$ which is invariant under $\f$. Again
this notion is independent of the choice of a negatively curved metric on $M$
and coincides with the notion used above.

Thus for a H\"older class $\zeta$ and a measure $\mu\in {\cal M}$
we can define $\int \zeta d\mu = \int \zeta(v,1)d\mu(v)$. Then the
{\sl pressure} $pr ([\zeta])$ of $[\zeta]$ is defined to be the supremum of
the values $h_\mu - \int
\zeta d\mu$ where $\mu$ ranges through $\m$ and $h_\mu$ is the {\sl entropy}
of
$\mu$.  If $[\zeta]$ is {\sl positive}, i.e. if $[\zeta]$ can be represented
by a cocycle $\bar\zeta$ which satisfies $\bar\zeta (v, t) > 0$ for all $v \in
T^1 M$ and all $t > 0$, then there is a unique number $a > 0$ such that $[a
\zeta]$ is {\sl normalized} i.e. that $pr ([a \zeta] ) = 0$.
The Gibbs current $[\mu_\zeta]$ of $[\zeta]$ then is the projection 
of the Gibbs equilibrium state of $a\zeta$ to ${\cal G}\tm$.

As was shown in [B2], 
for a surface $M$ there is a natural
bilinear form, the so called {\sl intersection form}, on the 
compact convex space
of geodesic currents
for $M$. This bilinear form is
continuous with respect to the
weak$^*$- topology on the space of geodesic currents.
The intersection of the Lebesgue Liouville current
of a negatively curved metric $g$ with any current $\beta$ is just
twice the integral of the length cocycle of $g$ with respect
to $\beta$. 

One
consequence of Theorem $B$ is the following: For every closed manifold
$M$, the map which assigns to
the class of the length cocycle of a strictly  $1/4$-pinched negatively curved
metric
on $M$ its Lebesgue Liouville current is injective. Thus we can
define the intersection between a current $\nu$ and the Lebesgue Liouville
current $\lambda_\ell$ of the length cocycle $\ell$ of such a metric
as
$\int \ell d\nu$.

We conjecture that the class of the length cocycle of every negatively curved
metric is determined by its Lebesgue-Liouville current, and that it is
possible to extend the above function to a continuous non-symmetric bilinear
form defined on the vector space spanned by all Gibbs currents.

%\hfill\eject
\bigskip
\bigskip
\centerline{\bf 2. Gibbs currents for surfaces}
\bigskip
\bigskip
To begin with, let $M$ be an arbitrary closed Riemannian manifold of negative
sectional curvature. Recall
that the space of  cohomology classes
of H\"older cocycles for the geodesic flow on $T^1M$
is independent of the choice of a metric of negative curvature on $M$
(compare [L] and the introduction).

We will only consider flip invariant H\"older cohomology classes. The
cohomology class $[f]$ defined by a function $f$ on $T^1M$ is called
{\sl positive} if it can be represented by a positive
H\"older function. This notion coincides with the one given
in the introduction.
 Denote by $\ch$ the space of positive flip invariant H\"older
classes. This space carries a natural (non-complete) topology as follows:

Recall that a {\sl geodesic current} for $M$ is a $\Gamma = \pi_1
(M)$-invariant locally finite Borel measure on
$\partial \tm \times \partial \tm - \Delta$ where $\Delta$ is the diagonal in
$\partial \tm \times \partial \tm$, which in addition is invariant under
the natural involution of
$\partial \tm \times \partial \tm$ exchanging the two factors. Geodesic
currents correspond naturally to finite $\Phi^t$-invariant Borel measures on
$T^1 M$ which
are invariant under the flip. We equip the space $\c$ of geodesic currents
with
the weak*-topology. With this topology, $\c$ is a locally compact space.

For every $\alpha \in \ch$ and every current $\eta \in \c$ the integral
$\int \alpha d\eta$ is well defined. We equip $\ch$ with the coarsest topology
such that the function $(\alpha, \eta) \in \ch \times \c \to \int \alpha d
\eta$ is continuous. This then corresponds to the topology of uniform
convergence for continuous functions $f$ on $T^1 M$ if we represent a current
by a $\Phi^t$-invariant finite Borel measure on $T^1 M$ and a
H\"older class by a H\"older continuous function on $T^1 M$.

Consider again the geodesic flow $\Phi^t$ on the unit tangent bundle $T^1 M$
of $M$. Let $\m$ be the compact convex space of $\Phi^t$-invariant
Borel-probability measures on $T^1M$ equipped with the weak*-topology. Recall
that the {\sl pressure} $pr (f)$ of a continuous function $f$ on $T^1 M$
is defined by
$pr (f) = \sup \{h_\nu - \int f d\nu \mid \nu \in \m\}$ where $h_\nu$ is the
entropy of $\nu$. Call a H\"older class {\sl normalized} if it can be
represented
by a function $f$ with $pr (f) = 0.$ Every normalized H\"older class is
positive. Vice versa, if $\alpha$ is a positive H\"older class, then there is a
unique constant $h (\alpha) > 0$ such that $h (\alpha) \alpha$ is normalized.
We call $h (\alpha)$  the {\sl topological entropy} of $\alpha$.

Identify the space of flip invariant normalized H\"older classes with the
projectivization $\p\ch$ of $\ch$.

Every element from $\p\ch$ determines its Gibbs current whose class
in the space $\p\c = \c /{{\bf R}_+}$ of {\sl projective currents}
where ${\bf R}_+$ acts on $\c$ via multiplication with a positive constant does
not depend on any choices made.
Thus we obtain a natural map of the space $\p\ch$ of projective flip
invariant positive H\"older
classes into the space $\p\c$ of projective currents.
We call the image of $\p\ch$ under this
map the {\sl space of projective Gibbs currents}. Also we call a current a
{\sl Gibbs current} if it corresponds to the Gibbs equilibrium state of a
normalized flip invariant H\"older class on $T^1 M$. Denote by $[\mu_\alpha]$
the class in $\p\c$ defined by the Gibbs state of $h (\alpha) \alpha$ where
$\alpha \in \ch$ and as before, $h (\alpha)$ is the topological entropy of
$\alpha$.

Clearly $h (a \alpha) = a^{-1} h (\alpha)$ for all $a > 0$,
moreover if $\alpha$ is the length cocycle of a negatively curved Riemannian
metric $g$ on $M$ then $h (\alpha)$ is just the topological entropy of the
geodesic flow of $g$. First of all we have:
\bigskip
{\bf Lemma 2.1:}
{\it $h (\alpha + \beta) \le h(\alpha)h(\beta)/(h (\alpha) + h (\beta))$
with equality if and only
if $\alpha$ is a constant multiple of $\beta$.}
\bigskip
{\it Proof:} Fix a negatively curved metric $g$ on $M$ with geodesic flow
$\Phi^t$.
Let $\alpha, \beta$ be normalized H\"older classes.
Let $\mu$ be the unique $\Phi^t$-invariant Borel probability measure on the
unit tangent bundle for $(M, g)$ which is a Gibbs equilibrium state for
the normalization of $\alpha + \beta$ and denote by
$h_\mu$ the entropy of $\mu$.
Then we have $\int \alpha d\mu \ge h_\mu/h(\alpha), \int \beta d\mu \ge
h_\mu/h(\beta)$ with equality if and only if $\alpha = \beta$.
Then $\int [ \alpha +  \beta]
d\mu \ge h_\mu (1/h (\alpha) + 1/h (\beta))$ and hence
the lemma follows from the fact that $h_\mu = h (\alpha + \beta) \int (\alpha
+ \beta) d\mu$.
{\bf q.e.d.}
\bigskip
We equip $\p\c$ with the topology induced from the weak*-topology on
$\c$.
\bigskip
{\bf Lemma 2.2:} {\it i) The topological entropy $h \colon \ch \to (0 ,
\infty)$ is continuous.
ii) The map $[\alpha] \in \p \ch \to [\mu_\alpha] \in \p \c$ is continuous.
\bigskip
Proof:} Let $\{\alpha_i\} \subset \ch$ be a sequence converging to some
$\alpha
\in \ch$. For $i \ge 1$ let $\eta_i \in \m$ be the unique Gibbs equilibrium
state for $h (\alpha_i) \alpha_i$. Since $\m$ is compact we may assume by
passing to a subsequence that the measures $\eta_i$ converge weakly to some
$\eta \in \m$. If $h_\nu$ denotes again the entropy of $\nu \in \m$ then $\nu
\to h_\nu$ is upper semi-continuous and hence $h_\eta \ge \lim_{i \to \infty}
\sup h_{\eta_i}$.

By the definition of the topology on $\ch$ we have
$$\int \alpha_i d\eta_i \to \int \alpha d\eta > 0 (i \to \infty),$$
in particular $\{h (\alpha_i)\}$ is bounded from above and below by a positive
constant. By passing to a subsequence we may assume that
$h (\alpha_i) \to \bar h > 0 (i \to \infty)$. Then
$h_\eta- \bar h\int \alpha d \eta \ge \limsup_{i \to \infty} h_{\eta_i}- h
(\alpha_i) \int \alpha_i d \eta_i = 0$ and consequently
$\bar h \le h (\alpha)$. On the other hand, if $\bar h < h (\alpha)$ then
there is
$\nu \in \m$ with $h_\nu - \bar h \int \alpha d\nu > 0$. Then also
$h_\nu - h (\alpha_i)\int \alpha_i  d\nu > 0$ for $i$ sufficiently large which
is impossible. In other
words, the function $h$ is indeed continuous. Moreover, if $\eta\in
{\cal M}$ is as
above, then $h_\eta - \bar h \int \alpha d \eta = 0$ and hence the current
determined by $\eta$ is contained in the class of $[\mu_\alpha]$. From this
the continuity of the map $[\alpha] \to [\mu_\alpha]$ is immediate.
{\bf q.e.d.}
\bigskip

Recall from [H5] that a {\sl cross-ratio} for $\Gamma = \pi_1 (M)$ is a
H\"older
continuous positive function $Cr$ on the space of quadruples of pairwise
distinct points in $\partial \tm$ with the following properties:

1) $Cr$ is invariant under the action of $\Gamma$ on $(\partial \tm)^4$.

2) $Cr (a, a', b, b') = Cr (a', a, b, b')^{-1}$

3) $Cr (a, a', b, b') = Cr (b, b', a, a')$

4) $Cr (a, a', b, b') Cr (a', a'', b, b')= Cr (a, a'', b, b')$

5) $Cr (a, a', b, b') Cr (a',b,a, b')Cr (b,a,a', b') = 1.$

Property 5) above is a consequence of properties 1) - 4) and the fact
that $Cr$ admits a H\"older continuous extension to the space of quadruples $(a,
a', b,b')$ of points in $\partial \tm$ which satisfy $\{a, a'\} \cap \{ b,
b'\} = \emptyset.$ This extension equals $1$ for every quadruple $(a, a', b,
b')$ for which either $a = a'$ or $b = b'$ (this was communicated to me by F.
Ledrappier).

We showed in [H5] that there is natural bijection between the space of
cohomology
classes of 
flip invariant H\"older cocyles and the space of cross ratios on $\partial \tm
\times
\partial \tm-\Delta$. We call a cross ratio $Cr$ {\it positive} if it
corresponds under this identification to a positive H\"older
 class.
\bigskip

{\bf Lemma 2.3:} {\it Let $Cr$ be a cross ratio on $\partial \tm$.
Then $Cr (a, b, c, d)^{-1} + Cr (b, c, d, a)^{-1}\hfil\break
\to 1$ as $a \to b$, locally uniformly in $(c, d)$ if and only if $Cr$ is
positive.}
\bigskip
{\it Proof:} Recall from [H5] that there is a H\"older continuous symmetric
function $(\;,\;)$ on $\partial \tm \times \partial \tm - \Delta$ (where
$\Delta$ denotes the diagonal) such that $\log Cr (a, b, c, d) = (a, d) +
(b, c) - (a, c) - (b, d)$ for pairwise distinct points $a, b, c, d$ in
$\partial \tm$. If $Cr$ is positive then $(a, b) \to \infty$ as $a \to b$ in
$\partial \tm$. Now
$$Cr (a, b, c , d)^{-1} + Cr (b, c, d, a)^{-1} = e^{(a, c)
+ (b, d)} [e^{- (a, d) - (b, c)} + e^{- (b, a) - (c, d)} ]$$
and this converges to $1$ as $a  \to b$ if and only if
$e^{- (b, a) - (c, d)} \to 0$ as $a \to b$. The
formula also shows that this convergence is locally uniform in $(c, d)$ if
$Cr$ is positive.
{\bf q.e.d.}
\bigskip
%Recall that a {\sl signed geodesic current} is a signed measure on
%$\partial \tm \times \partial \tm - \Delta$
%which can be written in the form $\mu_1-\mu_2$ where $\mu_1$ and $\mu_2$ are
%geodesic currents.

From now on we specialize to the case that $M$ is an oriented surface. The
ideal
boundary of its universal covering is naturally homeomorphic to $S^1$. Fix
once and for all an orientation for the circle $S^1$. This orientation then
determines for every ordered pair $(a, b)$ of points $a \ne b$ in $S^1$ a
unique half-open interval $[a, b[\subset S^1$ with endpoints $a$ and $b$.
Call a quadruple $(a_1, a_2, a_3, a_4)$ of pairwise distinct points in $S^1$
{\sl ordered} if $a_i$ is contained in $[a_{i - 1}, a_{i + 1} [$ for $i = 1,
\dots, 4.$

\bigskip
{\bf Proposition 2.4:} {\it For every Gibbs current $\mu$ there is a unique
cross ratio $Cr$ for $\Gamma = \pi_1 (M)$ such that $\log Cr (a, b, c, d) =
\mu [a, b[ \times [c, d[$ for every ordered quadruple $(a, b, c, d)$ in $S^1$.
Moreover $Cr$ is positive.}
\bigskip
{\it Proof:} Let $(a, b, c, d)$ be an ordered quadruple of pairwise distinct
points in $S^1$ and define $[a, b, c, d] = \mu [a, b[\times [c, d[$. Since
$\mu$ is a current and hence invariant under exchange of the intervals $[a,
b[$ and $[c, d[$ we have $[c, d, a, b] = \mu [c, d[ \times [a, b[ = [a, b, c,
d]$, i.e. $Cr = e^{[\;]}$ satisfies 3) above. Moreover, if $y \in ] a, b[$
then $(a, y, c, d)$ and $(y, b, c, d)$ are ordered and

$$\mu [a, b[ \times [c, d[ = \mu [a, y[ \times [c, d [ + \mu [y, b[ \times
[c, d[ \leqno (*)$$

which corresponds to $4)$ above for $Cr$ on ordered quadruples.

So far we have not used that $\mu$ is a Gibbs current; this now is essential
to show that [ ] is H\"older continuous. Fix a base point $x \in \tm$ and view
$S^1 = \partial \tm$ as the unit sphere in $T_x \tm$. By formula (*) it is
enough to show the following: For $(c, d) \in S^1 \times S^1 - \Delta$ and $b
\in ] d, c [$ there are constants $\beta > 0, \alpha > 0$ such that $\mu [a,
b[\times [c, d[ \le \beta \vert a- b\vert^\alpha$ whenever $a$ is sufficiently
close to $b$.

For this choose a positive H\"older function $f$ on $T^1 M$ such that $\mu$ is
the Gibbs current for $f$. Recall from [H5] that $f$ induces a symmetric
function $\alpha_f$ on $S^1 \times S^1 - \Delta$ such that $\beta^{-1}
\alpha_f (a, b)^{1 / \chi} \le (a \mid b) \le \beta \alpha_f (a, b)^\chi$
where $(\; \mid \;)$ is the Gromov product on $S^1 \times S^1 - \Delta$ with
respect to the base-point $x$ which
defines the H\"older structure on $S^1 = \partial \tm$ and $\chi \in (0, 1),
\beta > 0$ are fixed constants. Moreover for $a$ sufficiently close to $b, \mu
[a, b[\times[c, d[$ is bounded from above by a constant multiple of
$\alpha_f (a, b)$. From this H\"older continuity of $[ \; ]$ on ordered
quadruples of points in $\partial \tm$ is immediate.

Next, if $(a, b, c, d)$ is ordered, then we define $[b, a, c, d]$$ = - [a, b, c,
d]$; this then implies also that $[a, b, c, d] = - [a, b, d, c].$ Finally if
we put $[b, a, d, c] = [a, b, c, d]$ whenever $(a, b, c, d)$ is ordered then
we obtain an extension of [ ] to all quadruples of pairwise distinct points in
$\partial \tm$ which is independent of the choice of an orientation for
$\partial \tm$. Moreover $Cr = e^{[\; ]}$ clearly satisfies all defining
properties of a cross ratio.
{\bf q.e.d.}

\bigskip
We call the cross ratio $Cr$ defined as above by a Gibbs current $\mu$ the
{\sl intersection cross ratio} of $\mu$ and we write $[\; ]_\mu = \log Cr$. To
justify this notion, recall that the {\sl intersection form} $i$ is a
continuous bilinear form on the space of geodesic currents where the space of
currents is equipped with the weak*-topology (see [B2]).
Recall also that every free homotopy class $[\gamma]$ in $M$ defines a
unique geodesic current
which we denote again by $[\gamma]$ as follows: Represent $[\gamma]$ by a
closed
geodesic in $M$. The lifts of this geodesic to $\tm$ define a $\Gamma = \pi_1
(M)$-invariant subset of the space
$\partial \tm \times \partial \tm - \Delta$
of geodesics whose intersection with every compact subset of
$\partial \tm \times \partial \tm - \Delta$
is finite. Then $[\gamma]$ is the sum of all Dirac masses on all points of
$\partial \tm \times \partial \tm - \Delta$
corresponding to those lifts. Recall also that a free homotopy class in $M$ is
nothing else but a conjugacy class in $\pi_1 (M)$.

Then we have:
\bigskip
{\bf Lemma 2.5:} {\it Let $\mu$ be a Gibbs current, let $\Psi\in \Gamma$ and
denote by $a, b$ the fixed points of the action of $\Psi$ on $S^1$. Then for
every $\xi \in ] b, a[, $ the intersection of $\mu$ with the conjugacy class
of $\Psi$ equals
$$\vert [a, b, \xi, \Psi \xi]_\mu\vert.$$
}
\bigskip
{\it Proof:} By definition, if $\xi$ is such that $(a, b, \xi, \Psi \xi)$
is ordered, then
$$[a, b, \xi, \Psi \xi]_\mu = \mu [a, b [\times [ \xi, \Psi\xi [;$$
but the
right hand side of this equation just equals the intersection of
$\mu$ with current defined by the conjugacy class of $\Psi$ (see [O]).
{\bf q.e.d.}
\bigskip
Using Lemma 2.5, we observe next that the intersection cross ratio just
describes the duality between the Lebesgue-Liouville current of a metric $g$
of negative curvature on $M$ and the Bowen-Margulis current of the geodesic
flow for $g$. For this we recall again that there is a natural
1-1-correspondence between flip invariant H\"older classes and cross ratios ([H5
]).

\bigskip
{\bf Corollary 2.6:} {\it Let $\lambda_g$ be the Lebesgue Liouville current of
a
metric $g$ on $M$ of negative curvature. Then the intersection cross ratio of
$\lambda_g$ is twice the cross ratio defined by the length element of $g$.}
\bigskip
{\it Proof:} The computations in [B1] show that the intersection of the
Lebesgue Liouville current $\lambda_g$ of $g$ with a free homotopy class
$\gamma$ in $M$ is just twice the length of the closed geodesic with
respect to $g$ which represents $\gamma$. In other words, the value at
$\gamma$ of the H\"older cohomology class corresponding to the intersection
cross ratio of $\lambda_g$ equals twice the value at $\gamma$ of the
length cocycle of the metric $g$. Since a H\"older
class is determined by  its values on free homotopy classes, the corollary
follows.
{\bf q.e.d.}
\bigskip
Another consequence of Proposition 2.4 is the following.
\bigskip
{\bf Lemma 2.7:} {\it The map $[\alpha] \in \p \ch \to [\mu_\alpha] \in
\p\c$ is injective.
\bigskip
Proof:} Let $\alpha \in \ch$ be normalized and let $\mu \in [\mu_\alpha]$.
Let
$\Psi \in \Gamma = \pi_1 (M)$ and let $a$ be the attracting fix point for the
action of $\Psi$ on $\partial \tm,\; b$ be the repelling fix point. For fixed
$\xi \in ] b,a [$ we have $\alpha (\Psi, a) = 
-\lim_{k \to\infty} {1\over k}\log
[b, \xi, \Psi^k \xi,a]_\mu$ (see
[H5]) and hence $\alpha$ is completely determined by $[\;]_\mu$.
From this the lemma is immediate.
{\bf q.e.d.}

\bigskip
We can also ask for the reverse of the above procedure. Namely, let $Cr$
be a
cross ratio and write $[ \; ] = \log Cr$. For an ordered quadruple $(a, b, c,
d)$ of pairwise distinct points in $S^1$ define $\eta [a, b[ \times [c, d[ =
[a, b, c, d]$. By our assumption, this defines a function on those subsets of
$S^1 \times S^1 - \Delta$
which are products of non-empty
right-half open intervals in $S^1$. Since $[a, y, c, d] + [y, b, c, d] = [a,
b, c, d]$ for all $y \in ] a, b [$, this function has a natural extension to a
finitely additive function on the family of finite unions 
products of of right-half open
intervals in $S^1$. If this function is in addition positive, then it defines
a locally finite Borel measure on $S^1 \times S^1 - \Delta$
which we denote again by $\eta$. Since $[ \; ]$ is invariant under the action of
$\Gamma$ , the measure $\eta$ is in fact a geodesic current. We call
$\eta$ the {\sl intersection current} of $Cr$. In general however the finitely
additive function on finite unions of 
products of right-half open intervals in $S^1$ is not
$\sigma$-additive, which means that it is not a {\sl signed geodesic current}
(the formal difference of two geodesic currents). However we always call
$\eta$ the intersection current of $Cr$.

The results of [H5] show that every
H\"older class $\alpha$ on $T^1 M$ defines a unique intersection current
$\nu_\alpha$.
The assignment $\alpha \to \nu_\alpha$ is linear, i.e. if $\alpha,
\beta$ are H\"older classes and if $a, b \in {\bf R}$ then
$\nu_{a \alpha + b \beta} = a\nu_\alpha + b \nu_\beta.$

Clearly if the intersection current $\nu_\alpha$ is positive, i.e. if
$\nu_\alpha$ is in fact a current, then $\alpha$ is a positive H\"older class.
The reverse however is not true. To see this let $\alpha_1, \alpha_2$ be the
H\"older classes of the length cocyles of hyperbolic metrics $g_1 \ne g_2$ on
$M$. Since the metrics $g_1, g_2$ are bilipschitz equivalent there is a number
$\ep > 0$ such that $\alpha_1 - \ep \alpha_2$ is a positive H\"older class.

By Corollary 2.6 the intersection current $\lambda_i$ of $\alpha_i$ is ${1
\over 2}$ times the Lebesgue Liouville current of $g_i (i = 1 , 2)$. Since the
currents $\lambda_1$ and $\lambda_2$ are singular, the intersection current
$\lambda_1 - \ep \lambda_2$ of the positive H\"older class
$\alpha_1 - \ep \alpha_2$ is not positive.

Even if the positive H\"older class $\alpha$ is such that its intersection
current $\nu_\alpha$ is positive and ergodic under the action of $\Gamma =
\pi_1 (M)$ we do not know whether $\nu_\alpha$ is in fact a Gibbs current.

For every geodesic current $\eta$ and every H\"older class $\alpha$
the integral $\int \alpha d\eta$ is well defined.
More precisely, recall from [B2] that for two currents $\eta, \eta'$
the intersection $i (\eta, \eta')$ is defined.
Then we have:

\bigskip
{\bf Lemma 2.8:} {\it For every positive H\"older class $\alpha$ for which
$\nu_\alpha$ is positive and
every current $\beta$ we have $\int \alpha d \beta = i (\beta, \nu_\alpha).$}
\bigskip
{\it Proof:} Let $\gamma$ be a current defined by a conjugacy class in $\pi_1
(M).$ Let $\Psi \in \pi_1 (M)$ be contained in this conjugacy class and let
$a, b$ be the fixed points of the action of $\Psi$ on $S^1$. Then $i (\gamma,
\beta) = \beta ]a, b [ \times [\xi, \Psi\xi [$ for every $\xi \in S^1 - \{a,
b\}$ and every current $\beta$ and hence by the definition of the map $\alpha
\to \nu_\alpha$ (see [H5])
we have $i (\gamma, \nu_\alpha) = \int \alpha d\gamma$. Since finite sums of
weighted Dirac currents
of conjugacy classes in $\pi_1 (M)$ are dense in the space of all currents
equipped with the weak*-topology and since the intersection form is continuous
with respect to the weak*-topology ([B2]) the lemma follows.
{\bf q.e.d.}

\bigskip
Every Riemannian metric $g$ on $M$ of negative curvature defines
its Lebesgue-Liouville current $\lambda_g$.
\bigskip
{\bf Lemma 2.9:} {\it The set of Lebesgue Liouville currents  of negatively
curved metrics is dense in the set of Gibbs currents.}
\bigskip
{\it Proof:}
Since every Gibbs current can be approximated by a sequence of Dirac
currents for conjugacy classes in $\pi_1 (M)$ it suffices to show that the
closure
of the set of all Lebesgue Liouville currents of metrics of negative curvature
in the space $\c$ of geodesic currents equipped with the
weak*-topology contains all currents defined by conjugacy classes in
$\pi_1 (M)$.

For this recall first of all that a free homotopy class in $M$ which can be
represented by a simple closed curve defines the projective class of a
measured geodesic lamination.
On the other hand, projective measured laminations form the boundary in $\p\c$
of the projectivizations of the Lebesgue Liouville currents of hyperbolic
metrics on $M$ ([B1]). Thus we may restrict our attention to free homotopy
classes which can not be represented by simple closed curves.

Fix a hyperbolic metric $g$ on $M$ and let $\gamma$ be a closed geodesic in
$(M, g)$ with self intersections. By eventually changing the metric
$g$ we may assume that $\gamma$ has only double points.

For simplicity we
consider only the case that $\gamma$ has exactly one double point. The
construction given below is also valid in the general case. Let
$\gamma\colon [0, T] \to M$ be a parametrization by arc length such that
$\gamma (0) = \gamma (\tau) = \gamma (T)$ for some $\tau \in (0, T)$ (in other
words, $\gamma (0)$ is the double point).

For $\ep > 0$ let $Q (\ep)$ be the $\ep$-neighborhood of $\gamma$ in $M$. For
sufficiently small $\ep,\; Q (\ep) - \gamma [0 , \tau]$ has exactly two
components, one of which, say $Q_1,$ is homeomorphic to an annulus. Choose a
normal vector field $t \to N (t)$ along $\gamma$ which points on $(0, \tau)$
inside $Q_1$. Let $m > 0$ be a large number; then there are $\ep_1 (m), \ep_2
(m) \in (0, \ep)$ such that
$\exp\; \ep_1 (m) N (1/m) = \exp\; \ep_1 (m) N (\tau - 1/m) \in Q_1$ (where
exp is the exponential map of $(M, g)$ ) and\hfil\break
$\exp \;( -\ep_2 (m) N (1/m)) = \exp\; \ep_2 (m) N (\tau + 1/m)$). Observe
that
the absolute values of $\ep_1 (m)$ and $\ep_2 (m)$ only depend on $m$ and the
angle between $\gamma' (0)$ and $-\gamma '(\tau)$ at $\gamma (0)$.

Thus for each sufficiently large $m \ge 1$ we obtain a closed neighborhood $A
(m)$ of $\gamma$ in $M$ with the following properties:

%{\setindent{iii)}
{
\item{i)} $A (m) \supset A (m + 1)$ and $\cap_{m \ge 1} A (m) = \gamma$
\item{ii)} $A (m) = \cup^4_{i = 0} A_i (m)$ where the sets $A_i (m)$ are closed with
piecewise smooth boundary and pairwise disjoint interior.
\item{iii)} $A_0 (m)$ is a geodesic quadrangle in $M$ containing $\gamma (0)$ in its
interior, with vertices

}

$$\eqalign{\exp\;\ep_1 (m) N (1/m) =  x_1 (m), & \; \exp \; - \ep_2 (m) N
(1/m) = x_2 (m),\cr \exp\;( -\ep_1 (m) N (T-1/m))
= x_3  (m) \hbox{ and } & \; \exp\; \ep_2 (m) N (T - 1/m) = x_4  (m).}$$

%{\setindent{iv)}
{
\item{iv)} The boundary of $A_i (m) (i = 1, \dots , 4)$ contains two smooth geodesic
segments of length $\ep_1 (m)$ or $\ep_2 (m)$ which meet at $x_i (m)$ and are
subarcs of the boundary of $A_0 (m)$. It also contains a subarc of
$\gamma$. With
respect to normal exponential coordinates based at $\gamma$ the metric on $A_i
(m)$ can be written in the form
$(\cosh s)^2 dt^2 + ds^2 (s \in [0, \ep_1(m) ]$ or $s \in [0, \ep_2 (m) ]).$

}

Change now the metric in the interior of $A (m)$ as follows: Fix $m \ge 1$
sufficiently large and for $i = 1, 2$ choose a diffeomorphism
$\Psi_{i, m} \colon [0, m] \to [0, \ep_i (m) ]$
with the following properties:
\hfill\eject
%{\setindent{2)}
{
\item{1)} $\Psi_{i, m}' (s) = 1$ for $s$ near $0$ and $s$ near $m$.
\item{2)} ${-\Psi_{i, m}'' (s) \over\Psi_{i, m}' (s)} <
{\cosh (\Psi_{i, m} (s))\over \sinh (\Psi_{i, m} (s))}$
for $s \in [0, m]$

}

Property 2) can be fulfilled since
${\cosh  (s)\over \sinh (s)}\cdot s \to 1$ as $s \searrow 0.$
Define now a new metric $g_m$ on $M $ as follows:

%{\setindent{b)}
{
\item{a)} $g_m$ coincides with the original metric $g$ outside $A (m)$.
\item{b)} In normal exponential coordinates based at $\gamma$ the restriction
of $g_m$ to $A_j (m)$ can be written in the form $dg_m = (\cosh s)^2 dt^2 +
[(\Psi_{i, m}^{-1})' (s)]^2 ds^2 (i = 1 {\rm \;or}\; i = 2\; {\rm and\;}
j = 1, \dots , 4).$

}

This defines $g_m$ on $M - A_0 (m)$. Observe that $g_m$ has negative curvature
(this is guaranteed by property 2) for $\Psi_{i, m}$ ) and has a natural
extension to the boundary of $A_0$ which coincides with $g$ near the vertices
$x_j (m) (j = 1, \dots, 4)$. With respect to this extension,
$A_0$ is a geodesic quadrangle of side length $2 m$ and such that the sum of
the
internal angles of the quadrangle is strictly less than $2\pi$. This means
that
we can extend $g_m$ to a metric of negative curvature on $A_0 (m)$ in such a
way that as $m \to \infty$ the geodesic quadrangle $(A_0 (m), g_m)$ approaches
an euclidean parallelogram of side length $2 m$ whose angles are determined by
the angle between $\gamma' (0)$ and $- \gamma' (\tau).$ Moreover we may assume
that $\gamma$ remains a $g_m$-geodesic for all $m$.

Let $\lambda_m'$ be the Lebesgue Liouville current of $g_m$ and write
$\lambda_m = {1 \over 4 m} \lambda_m'$. Choose a free homotopy class in $M$
which is prime, different from the class of $\gamma$ and which is represented
by the $g$-geodesic $\eta$.
 By abuse of notation we use the same symbol
for
$\eta$ and its free homotopy class.

After a slight perturbation of $\eta$ we may
assume that $\eta$ does not pass
through $\gamma (0).$ Then $\eta$ intersects $\gamma$ in finitely many points
$\eta (s_1), \dots, \eta (s_q)$ where $q = i (\gamma, \eta)$. For sufficiently
large $m > 1$ the intersection of $\eta$ with $A (m)$ has then precisely $q$
connected components, and the $g_m$-length of each of these components is
contained in the interval $[2m, 2m + c]$ where $c > 0$ is a constant which
depends on $\eta$ but not on $m$.
Thus $i (\lambda_m, \eta) = {1 \over 4m} i (\lambda_m', \eta) \le {1 \over 2m}
\cdot g_m - {\;\rm length \; of\;} \eta \to q (m \to \infty).$

In other words, we have $\limsup_{m \to \infty} i (\lambda_m, \eta) \le i
(\gamma, \eta)$.

Since
the intersection form $i$ on $\c \times \c$ is continuous and since the
set of finite weighted sums of Dirac
currents of free homotopy classes different from $\gamma$ is dense in $\c$
we conclude that the sequence
$\{\lambda_m\}$
is bounded in $\c$. By passing
to a subsequence we may therefore assume that
$\{\lambda_m\}$ converges weakly to a current $\lambda \in \c$. Then
$i (\lambda, \eta) \le i (\gamma, \eta)$
for every free homotopy class $\eta \not= \gamma$
and therefore $\lambda \le \gamma$. On the
other hand the current $\gamma$ is ergodic and hence $\lambda = a \gamma$ for
some $a \ge 0$.

We are left with showing that $a \ne  0$. For this choose again a free
homotopy class $\eta$ with
$i (\gamma, \eta) = q > 0$. Then every curve in $M$ representing $\eta$
intersects $\gamma$ in at least $q$ points, where an intersection at the point
$\gamma (0)$ has to be counted twice. But this means that the $g_m$-length of
every such curve is at least $2m$ and hence $\lim\inf_{m \to \infty} i
(\lambda_m, \eta) = i (\lambda, \eta) > 0$. This shows $a \ne 0$ and finishes
the proof of the lemma.
{\bf q.e.d.}
\bigskip

To summarize the considerations in the beginning of this chapter we have the
following situation: To every element $[\alpha]$ in the space $\p\ch$ of
projective flip invariant positive H\"older classes we can associate the
projective class $[\mu_\alpha]$ of the Gibbs equilibrium state defined by
$\alpha$ and also the projective class $[\nu_\alpha]$ of the intersection
current of $\alpha$. This defines two injective continuous maps $[\alpha] \to
[\nu_\alpha]$ and $[\alpha] \to [\mu_\alpha]$ of $\p\ch$ into the space of
projective classes of finitely additive signed measures on ${\cal G} \tm$.
The image of
$\p\ch$ under the second map is just the space of projective Gibbs currents.

With this terminology, the result of Katok in [K ] can be formulated as
follows:
If $\alpha$ is the class of the length cocycle of a metric on $M$ of negative
curvature, then $[\mu_\alpha] = [\nu_\alpha]$ if and only if this metric is of
constant negative curvature.

Recall the definition of the {\sl topological entropy} $h (\alpha)$ of a
positive H\"older class $\alpha$.

The proof of Katok uses another result proved in the same paper which can be
slightly generalized as follows
(for a generalization to higher dimensions see [B-C-G]):

\bigskip
{\bf Theorem 2.10:} {\it Let $M$ be a compact surface of Euler characteristic
$\chi (M) < 0$ and let $\alpha$ be a positive H\"older class such that
$\nu_\alpha$ is a Gibbs current with $i
(\nu_\alpha, \nu_\alpha) = \pi^2 \vert \chi (M) \vert$.
Then $h (\alpha) \ge {1 \over 4}$ with equality if and only if $\alpha$
is the length cocycle of a hyperbolic metric on $M$.}
\bigskip
{\it Proof:} Let $g, \tg$ be conformally equivalent Riemannian metrics on $M$
of negative curvature, i.e. there is a function $a \colon M \to (0, \infty)$
such that $a g = \tg$. Assume that the $g$-volume of $M$ and the $\tg$-volume
equals $(2\pi)^{-1}$. If $\lambda_g, \lambda_{\tg}$ denotes the Lebesgue
Liouville current of $g, \tg$ then this just means that $i (\lambda_g,
\lambda_g) = 2 = i (\lambda_{\tg}, \lambda_{\tg})$ (see [B1]).

Let $V$ be the unit tangent bundle of $(M, g)$ and let $\lambda$ be the
Lebesgue Liouville measure of $g$ on $V$. For a $g$-unit vector $v \in V$
denote by $\rho (v)$ the $\tg$-norm of $v$. If $P \colon V \to M$ denotes the
canonical projection then $\rho (v) = (a (Pv))^{1/2}$ and consequently

$$\int_V \rho d \lambda = \int_V (a^{1/2} \circ P) d\lambda \le (\int (a
\circ P) d\lambda)^{1/2} = 1$$
with equality if and only if $a \equiv 1$.

Let now $\alpha (g), \alpha (\tg)$ be the length cocycles of the metrics $g,
\tg$. Since clearly\hfil\break
$\int \alpha (\tg)\; d\lambda_g \le \int_V \rho
\; d\lambda$, Lemma 2.8 shows that
$$i (\lambda_g, \nu_{\alpha (\tg)} ) \le 1,\; i  (\lambda_\tg, \nu_{\alpha
(g)}) \le 1 \leqno(*)$$
with equality if and only if $g = \tg$.

Assume now in addition that $\tg$ is a metric of constant curvature and let $h
(g), h(\tg)$ be the topological entropy of the geodesic flow for $g, \tg$.
Then
the cocycles $\alpha (g) h (g)$ and $\alpha (\tg) h (\tg)$ are normalized.
Now $\lambda_{\tg} = 2 \nu_{\alpha (\tg)},$ moreover the Lebesgue Liouville
measure of
$\tg$
equals the Gibbs equilibrium state for $\alpha (\tg) h (\tg)$ and therefore

$$\eqalign{\int\alpha (g) h (g)\; d \lambda_\tg = h (g)i(\lambda_\tg,
\nu_{\alpha (g)}) \ge\cr
\int h (\tg)\alpha(\tg)\; d \lambda_\tg = h (\tg)i(\lambda_\tg, \nu_{\alpha
(\tg)})=\cr
h (\tg){1 \over 2}i(\lambda_\tg, \lambda_\tg) = h (\tg).\cr}\leqno(**)$$

Together with (*) this shows that $h (g) \ge h (\tg)$ with equality if and
only if
$g
=
\tg$.

Notice that $h (\tg)$ does not depend on the metric $\tg$ of constant
curvature with \hfil\break
$i (\lambda_{\tg}, \lambda_{\tg}) = 2$ (compare [B1]). We write $h > 0$ for
this common number.
Let now $\alpha$ be a positive H\"older class such that $\nu = \nu_\alpha$ is a
Gibbs
current and that $i (\nu_\alpha, \nu_\alpha) = {1 \over 2}$. By Lemma 2.9,
Lebesgue Liouville currents of negatively curved metrics are dense in the set
of Gibbs currents. This means that there is a sequence $\{g_j\}_j$ of
negatively curved metrics on $M$ with Liouville currents
$\lambda_{g_j}$ and
such that $\lambda_{g_j} \to 2\nu$. Since the intersection form is continuous
we may assume that $g_j$ is normalized in such a way that
$i( \lambda_{g_j}, \lambda_{g_j}) = 2.$

Let $\alpha_j$ be the H\"older class determined by the length cocycle of the
metric $g_j$.
Since the intersection form $i$ is continuous on the space of geodesic
currents we obtain from Lemma 2.8 and the fact that 
$\nu_{\alpha_j}\to \nu_\alpha$ weakly in ${\cal C}$ that $\alpha_j\to
\alpha$ in ${\cal H}$. Lemma 2.2 then shows that 
$h(\alpha_j)\to h(\alpha)$ and therefore from
the above consideration we conclude that $h (\alpha) \ge h (\tg) =
h$ where
$\tg$ is a metric  of constant curvature with $i (\lambda_\tg,
\lambda_\tg)=2$.

Assume now that $h (\alpha) = h$.
Let $\tg_j$ be a metric of constant 
curvature on $M$ which is conformally equivalent to
$g_j$ and such that $i (\lambda_{\tg_j}, \lambda_{\tg_j}) = 2.$
By our assumptions we have $h (g_j) - h(\tg_j) \to 0 (j \to \infty).$
By (*) and (**) above, this means that
$i (\lambda_{g_j}, \nu_{\alpha(\tg_j)} ) \to 1 (j \to \infty).$

Recall that the space of projective currents $\p\c$ is compact; hence by
passing to a subsequence we may assume that the projective classes
$[\lambda_{\tg_j}]$
of
$\lambda_{\tg_j}$ converge weakly to a projective current $[\beta]$. Then
$\{{\tg_j}\}$ is an unbounded sequence in Teichm\"uller space if and only if
$\beta$ is a geodesic lamination, i.e. it satisfies $i (\beta, \beta) = 0$
([B1]).

On the other hand, $\lambda_{g_j} \to \nu_\alpha$ and $\nu_\alpha$ is a Gibbs
current, i.e. its intersection with every current is positive. If $b_j > 0$ is
such that $b_j \lambda_{\tg_j} \to  \beta \ne 0$ then
$i (\lambda_{g_j}, b_j \lambda_{\tg_j}) =
b_j i (\lambda_{g_j}, \lambda_{\tg_j}) \to i (\nu_\alpha, \beta) > 0 \,
(j \to \infty)$, again by continuity of the intersection form. But
$i (\lambda_{g_j},\lambda_{\tg_j}) \to 2 (j \to \infty)$
and consequently $\{ b_j \}$ is bounded from below by a positive number. This
implies in turn that $\{ \tg_j \}$ is bounded in Teichm\"uller space and hence
that
$[\lambda_{\tg_j}]\to[\lambda_{\tg}]$
for some metric $\tg$ of constant curvature.

If we normalize $\lambda_{\tg}$ in such a way that
$i (\lambda_{\tg}, \lambda_{\tg}) = 2$ then
$\lambda_{\tg_j}\to \lambda_{\tg}$ and also
$\lambda_{g_j}\to \lambda_{\tg}$ and consequently
$\nu_\alpha = {1\over 2} \lambda_{\tg}$. This finishes the proof of the
theorem.
{\bf q.e.d.}
\bigskip
If $\alpha$ is a positive H\"older class such that $\nu_\alpha$ is a Gibbs
current, then we can also define the {\sl metric entropy} $h_m (\alpha)$ of
$\alpha$ as follows: Let $\beta$ be the unique normalized positive H\"older
class such that
$[\mu_\beta] = [\nu_\alpha]$ and define
$h_m (\alpha) = \int
\beta d \nu_\alpha/
i (\nu_\alpha , \nu_\alpha)=
i (\nu_\alpha , \nu_\beta)  /
i (\nu_\alpha , \nu_\alpha)$.

If $\alpha$ is the length cocycle of a metric $g$ on $M$ of negative
curvature, then
$h_m (\alpha)$ equals the metric entropy of the geodesic flow
of $g$. Moreover
$h_m (a \alpha)= a^{-1}h_m ( \alpha)$ and if
$[\nu_\alpha] = [\mu_\alpha]$ then
$h( \alpha)= h_m ( \alpha)$.

In analogy to Theorem 2.10 we can also estimate
$h_m ( \alpha)$:

\bigskip
{\bf Theorem 2.11:} {\it Let $\alpha$ be a positive H\"older class for the
surface $M$ such that
\hfil\break  $i (\nu_\alpha , \nu_\alpha) = \pi^2 \vert \chi
(M)\vert$ and that
$\nu_\alpha$ is a Gibbs current. Then
$h_m ( \alpha) \le {1 \over 4}$ with equality if and only if $\alpha$ is the
length cocycle of a hyperbolic metric on $M$.}
\bigskip
{\it Proof:} As in the proof of Theorem 2.10, let $g, \tg$ be conformally
equivalent Riemannian metrics with Lebesgue Liouville currents
$\lambda_g, \lambda_{\tg}$ and length cocycles $\alpha (g), \alpha (\tg)$ and
such that
$i (\lambda_g, \lambda_g) = i (\lambda_{\tg}, \lambda_{\tg}) = 2, \tg$ of
constant curvature. Since the cocycle $
h (\alpha (\tg)) \alpha (\tg)$ is
normalized we have $\int h (\alpha (\tg)) \alpha (\tg) d \lambda_g =
 h (\alpha (\tg)) i (\lambda_g, \nu_{\alpha (\tg)})   \ge
 h_m (\alpha (g))\cdot {1 \over 2} i (\lambda_g, \lambda_g) =
 h_m (\alpha (g)).$ But $ i (\lambda_g, \nu_{\alpha (\tg)}) \le 1$
with equality only if $g = \tg$ by the proof of Theorem 2.10. Thus we obtain
the statement of the theorem
in the case that $\alpha$ is the class of a length cocycle for a
negatively curved Riemannian metric. The general case now follows as
in the proof of Theorem 2.10. 

Namely, let $\alpha$ be such that $\nu_\alpha$ is a Gibbs current
and $i(\nu_\alpha, \nu_\alpha)= {1\over 2}$. Let $\alpha_j$ be the
H\"older class of the length cocycle of a metric $g_j$ on $M$
such that $i(\nu_{\alpha_j}, \nu_{\alpha_j}) ={1\over 2}$ and
$\nu_{\alpha_j}\to \nu_\alpha$. Choose a metric ${\tilde g_j}$ on $M$
which is conformally equivalent to $g_j$ and satisfies $i(\lambda_
{\tilde g_j}, \lambda_{\tilde g_j})=2$; then
$i(\nu_{\alpha_j}, \lambda_{\tilde g_j})\leq 1$ for all $j$.

Now if $\tilde g_j$ is unbounded in Teichm\"uller space, then
there is a sequence $\{b_j\}\subset [0, \infty)$ such that
$b_j \to 0$ and that $\lambda_{\tilde g_j} b_j$ converges weakly to a
measured lamination $\beta$. Then 
$i(\nu_{\alpha_j}, b_j\lambda _{\tilde g_j}) = b_j
i(\nu_{\alpha_j}, \lambda_{\tilde g_j}) \to i(\nu_\alpha, \beta)>0$
which is impossible. Thus the sequence $\{\lambda_{\tilde g_j}\}$ is
bounded and by passing to a subsequence we may assume that 
$\lambda_{\tilde g_j}\to \lambda_{\tilde g}$ $(j\to \infty)$
where $i(\lambda_{\tilde g},\lambda_{\tilde g}) =2$ and $\tilde g$
is a metric of constant curvature. Then $i(\nu_{\alpha_j}, \lambda
_{\tilde g_j})\to i(\nu_\alpha, \lambda_{\tilde g})\leq 1$ and
hence $h_m(\alpha) \leq h(\alpha({\tilde g}))$ by the above 
argument. The discussion of equality is exactly the same as
the discussion of the equality case in the proof of Theorem 2.10.
{\bf q.e.d.}
\bigskip

Now if again $\alpha$ is a positive H\"older class such that
$[\nu_\alpha] = [\mu_\alpha]$, then in particular
$\nu_\alpha$ is a Gibbs current and hence Theorem 2.10 and 2.11 combined show:

\bigskip
{\bf Corollary 2.12:}
{\it $[\mu_\alpha] = [\nu_\alpha]$ for a positive H\"older class $\alpha$ if and
only
if $\alpha$ is the class of the length cocycle of a metric on $M$ of constant
negative curvature.}
\bigskip
{\it Remark:} In his paper [K], Katok showed Corollary 2.12 for length
cocycles of negatively curved Riemannian metrics. The arguments given here are
more formal and technically much easier, but follow the same principal idea
as the argument of Katok.
If $M$ is a compact rank 1 locally symmetric space, then the
analogue to Theorem 2.10 for length cocycles of negatively curved metrics is
due to Besson, Courtois and Gallot ([B-C-G]). However our simple estimate
using conformal equivalence of Riemannian metrics on $M$ is not valid any
more and the generalization in our Theorem 2.10 for arbitrary H\"older classes
is not true in higher dimensions (compare Section 3). Moreover, the analogue
of Theorem 2.11 fails even in the case that
the dimension of $M$ equals 3, as was pointed out by Flaminio ([F]). This is
however not so surprising given the fact that with the notation of the proof
of Theorem 2.11 we have $\int_V \rho d\lambda > \int_V d\lambda$ whenever $g
\ne \tg$ and $g$ and $\tg$ are metrics with the same volume form. However the
analogue of Corollary 2.12 may well be true for manifolds which carry a
locally symmetric metric.
\bigskip

Recall the following result of Bonahon ([B1]): If $\mu$ is a geodesic current
on a surface of negative Euler characteristic such that
$$e^{- \mu [a,b] \times [c,d]} + e^{-\mu [b, c] \times [d, a]} = 1$$
for every ordered quadruple $(a, b, c, d)$ then $\mu$ is the Lebesgue
Liouville current of a hyperbolic metric. This gives a purely algebraic
description of cross ratios "of maximal entropy". We can use the above
considerations to give a slightly sharper version of this result.

\bigskip

{\bf Corollary 2.13:} {\it Let $\mu$ be a geodesic current such that there is
a number $\rho > 1$ with $\rho^{-1}
e^{- \mu [a,b] \times [c,d]} \le
1 - e^{-\mu [b, c] \times [d, a]} \le \rho e^{- \mu [a,b] \times [c,d]} $ for
every ordered quadruple $(a, b, c, d)$ in $\partial \tm$. Then $\mu$ is the
Lebesgue-Liouville current of a metric of constant negative curvature.}
\bigskip
{\it Proof:} Let $\mu$ be a geodesic current which satisfies the assumptions
in the corollary. For two points $a \ne b \in S^1$ we then can define a
measure $\mu_{a b}$ on $] b, a [\subset S^1$ by $\mu_{ab} [x, y[ = \mu[a,b[
\times
[x, y[$.  We want to show that the measure class of $\mu_{ab}$ does not
depend on $a \ne b$.

For this choose $c \in ] a, b[$, let $x \in ] b, a [$ and let $y \in ] b, a [$
be a point near $x$. By our assumption there is a neighborhood $U$ of $x$ in
$] b, a [$ such that

$$1 - e^{-\mu_{cb} [x, y[}= 1 - e^{-\mu[c,b[\times [x, y[}\ge
\rho^{-1} e^{-\mu[b,x [\times[ y,c[}$$
for all $y \in U - \{x\}$ and moreover

$$e^{-\mu [b,x [\times [ y,c[}= e^{-\mu [b,x [\times [y,a[}
e^{-\mu [b,x [\times[a, c[}\ge
\rho^{-1}(1 - e^{-\mu[x,y [\times[ a,b[}) e^{-\mu [b,x [\times [a, c[}.$$

For $y$ sufficiently close to $x$ we have

$$\eqalign{1 - e^{-\mu_{cb} [x, y[}= &\mu_{cb} [x, y[ + o (\mu_{cb} [x,
y[),\cr
1 - e^{-\mu_{ab} [x, y[}= &\mu_{ab} [x, y[ + o (\mu_{ab} [x, y[ )}$$

and consequently the measures $\mu_{cb}$ and $\mu_{ab}$ are absolutely
continuous at $x$ with Radon Nikodym derivative

$$1 \ge {d\mu_{cb}\over d\mu_{ab}} (x) \ge \rho^{-2} e^{- \mu [b,x[ \times
[a,c[}.$$

This means that the current $\mu$, viewed as a $\Phi^t$-invariant Borel
measure on the unit tangent bundle $T^1 M$ of $M$, is absolutely continuous
with respect to the stable foliation, with locally bounded Radon Nikodym
derivative.

Let $\alpha$ be the H\"older class determined by $\mu$. Since $\mu$ is a
current, $\alpha$ is positive and is determined by its values on periodic
points for the geodesic flow on $T^1 M$. More precisely, for every $\Psi \in
\pi_1 (M)$ with attracting fix point $a$, repelling fix point $b$ for its
action on $S^1$ and every $\xi \in ] b,a [$ and $k \ge 1$ we have $\alpha
(\Psi, a) = {1 \over k} \mu [a, b[ \times [\xi, \Psi^k \xi [$ (where we view
$\alpha$ as a function on $\pi_1 (M)\times \partial \tm$).

On the other hand, $\lim_{k \to \infty} {1 \over k} \log \mu [b, \xi [\times [
\Psi^k \xi, a[$ is just the asymptotic logarithmic decay of the conditional
measure $\mu_{b \xi}$ at the attracting fix point $a$ under the action of
$\Psi$. For sufficiently large $k \ge 1$ we conclude that
$$\eqalign{ \rho^{-1} e^{-k \alpha(\Psi, a)} &\le 1- e^{-\mu [b,\xi [\times
[\Psi^k \xi,a [} =\cr
\mu [b,\xi [\times [\Psi^k \xi,a [ & + o (\mu [b,\xi [\times [\Psi^k \xi,a
[)\le
\rho e^{- k \alpha (\Psi, a)} }$$
and therefore
$-\alpha (\Psi,a ) = \lim_{k \to \infty} {1 \over k} \log \mu [b, \xi [\times
[ \Psi^k \xi, a[.$
But this just means that $[\mu] = [\nu_\alpha] = [\mu_\alpha]$ and hence the
corollary follows from Corollary 2.12.
{\bf q.e.d.}

\bigskip
\bigskip
\centerline{\bf 3. Contact cocycles and Lebesgue measures}
\centerline{\bf in higher dimensions}
\bigskip
\bigskip
In Section 2 we looked at signed geodesic currents for a closed hyperbolic
surface defined by the cross ratio of a H\"older class. This assignment can be
considered
as associating to the (de Rham) cohomology class of a H\"older continuous 1-form
on the unit tangent bundle of this surface
its derivative, viewed as a volume form on the
space of geodesics.

In this section we discuss an analogue of this in higher dimensions: Namely
the Lebesgue-Liouville current of a negatively
curved metric $g$ on a closed manifold $M$ can always be obtained in a purely
conbinatorical way
from the length cocycle of $g$. As a corollary we obtain Theorem B from the
introduction. For this we exploit the fact that the cross ratio of the length
cocycle can be viewed as a coarse version of the symplectic structure on the
space of geodesics.

To begin with, consider the vector space
${\bf R}^{2n}= {\bf R}^n \times {\bf R}^n$,
equipped with some smooth symplectic form
$\rho$. Assume that
for every fixed $x \in {\bf R}^n$ the sets $\{x\} \times {\bf R}^n$ and ${\bf R}^n
\times \{x\}$ are Lagrangian submanifolds of $({\bf R}^{2n}, \rho)$. We define a
function $\v$ on ${\bf R}^{2n}= {\bf R}^n \times {\bf R}^n$
as follows: For $x, y \in {\bf R}^n$ let $c_x \colon [0, 1] \to {\bf R}^n$ be a smooth
curve joining $c_x (0) = 0$ to $c_x (1) = x$ and choose similarly a curve $c_y
\colon [0, 1] \to {\bf R}^n$ joining $0 $ to $y$. Define a map $\Psi \colon [0,
1]^2 \to {\bf R}^{2n}$ by $\Psi (s, t) = (c_x (s), c_y (t))$ and let $\v (x, y) =
\int_{\Psi [0, 1]^2} \rho.$

\bigskip
{\bf Lemma 3.1:} {\it The function $\v$ does not depend on the choice of the
curves $c_x, c_y$.}
\bigskip
{\it Proof:} Choose a 1-form $\zeta$ such that $d \zeta = \rho$. Let $c_x,
\tc_x$ and $c_y, \tc_y$ be curves joining $0$ to $x$ and $y$ and let $\Psi,
\tp \colon [0, 1]^2 \to {\bf R}^{2n}$ be the corresponding maps defined as above.

By Stoke's theorem we have
$\int_{\Psi [0,1]^2} \rho = \int_{\partial\Psi[0, 1]^2} \zeta$
and
$\int_{\tp [0,1]^2} \rho = \int_{\partial\tp[0, 1]^2} \zeta$.
Now the oriented boundary $\partial\Psi - \partial \tp$ consists of 4 loops
(one of them is $c_x \circ \tc_x^{-1}$) contained entirely in a Lagrangian
submanifold for $\rho$. Another application of Stokes's theorem then shows
$$\int_{\partial\tp [0,1]^2} \zeta = \int_{\partial\Psi[0, 1]^2} \zeta.$$
{\bf q.e.d.}
\bigskip
If $\rho$ is the standard symplectic form $\rho_0 = \sum dx_i \wedge dy_i$ on
${\bf R}^n \times {\bf R}^n$, then we write $\v_0$ for the function as in Lemma 3.1.

Consider now ${\bf R}^{2n} = {\bf R}^n \times {\bf R}^n$ with the standard form $\rho_0$
and
the corresponding function $\v_0$. For a product set $Q_1 \times Q_2\subset
{\bf R}^n \times {\bf R}^n$ define the {\sl symplectic diameter}
$\delta (Q_1 \times Q_2)$ by $\delta (Q_1 \times Q_2) = \sup \{ \vert \v_0 (x,
y) \vert \mid x \in Q_1, y \in Q_2\}$ and for a compact subset $K$ of ${\bf R}^n$
and $r \ge 0$ let $B (K, r) = \{y \in {\bf R}^n\mid \vert \v_0 (x, y) \vert \le r$
for all $x \in K \}$. Clearly we have
$B (K, r) = r B (K, 1)$, moreover $B (K, r)$ is star-shaped about $0$ and
invariant under reflection at the origin.
Moreover we have:
\bigskip
{\bf Lemma 3.2:} {\it Let $K \subset {\bf R}^n$ be a compact neighborhood of $0$
with dense interior. Then $B (K, 1)$ is a compact convex body in ${\bf R}^n$ which
is reflection-symmetric at the origin. For every $x \in {\bf R}^n$ contained in
the boundary of $B (K , 1)$ a hyperplane $H \subset {\bf R}^n$ through $x$ is
supporting for $B (K,1)$ if and only if there is
$z \in K$ with $\v_0 (z, H - x) = 0$ and such that
$\vert \v_0 (z, x) \vert = \max \{ \vert \v_0 (y, x) \vert \mid y \in K \}$.
\bigskip
\it Proof:} For $x \in {\bf R}^n$ write $\alpha (x) = \sup \{ \vert \v_0 (z, x)
\vert \mid z \in K\}$. Since by assumption $K$ is a compact neighborhood of
$0$ in ${\bf R}^n$ we have $\alpha (x) > 0$ for $x \ne 0$ and hence $B (K, 1)$ is
compact with non-empty interior. Moreover $B (K, 1)$ is clearly star-shaped
about $0$.

Let $x \in B(K, 1)$ be a point from the boundary of $B (K, 1)$. Since $K$ is
compact there is $z \in K$ such that $\vert\v_0 (z, x) \vert = \alpha (x)
= 1$. Define $H (x) = \{ y \in {\bf R}^n \mid \v_0 (z, y) = 0 \}$. Then $H (x)$
is a hyperplane in ${\bf R}^n$ and $B (K, 1) \subset W (x) = \{ tx + y \mid - 1
\le t
\le 1,\  y \in H (x)\}.$

But this just means that $H (x)$ is a supporting hyperplane for $B (K, 1)$ at
$x$. Since $W (x)$ is convex and $B (K, 1) = \bigcap_{x \in \partial B (K, 1)}
W (x)$ we conclude that $B (K, 1)$ is convex.

For a dense set of points in the boundary of $B (K, 1)$ a supporting
hyperplane is unique. Therefore by continuity we conclude that the set of
supporting hyperplanes for $B (K, 1)$ at a point $x \in \partial B (K, 1)$ is
in $1-1-$correspondence with pairs of points $\pm z$ such that $z \in K$ and
that $\vert \v_0 (z, x) \vert = \alpha (x)$. This shows the lemma.
{\bf q.e.d.}
\bigskip
Let now $K \subset {\bf R}^n$ be a compact set with non-empty interior and write
$B_1 = B(K, 1),$
\hfil\break $ B_2 = B (B (K, 1), 1)$. By definition, $K \subset B_2$ and
$\vert \v_0 (x, y) \vert \le 1$ for all $x \in B_1, y \in B_2$. Moreover $B_1$
and $B_2$ are compact convex bodies in ${\bf R}^n$, reflection symmetric at the
origin. Now the symplectic form $\rho_0$ defines a linear isomorphism of ${\bf R}
={\bf R}^n \times \{0\}$ onto the dual of ${\bf R}^{n} = \{0\} \times {\bf R}^n$, and with
respect to this identification $B_2 \subset ({\bf R}^n)^*$ is just the dual of
the convex body $B_1 \subset {\bf R}^n$. We call $(B_1, B_2) \subset {\bf R}^n \times
{\bf R}^n$ a {\sl dual pair of convex bodies} in ${\bf R}^n$. If $B_1$ is an
ellipsoid, i.e. if $B_1$ is the image of the unit ball in ${\bf R}^n$ (equipped
with the standard inner product) under a linear isomorphism of ${\bf R}^n$, then
we call $(B_1, B_2)$ an {\sl euclidean dual pair}.

Now recall that the n-th exterior power $\rho^n_0$ of $\rho_0$ is a volume
form
on ${\bf R}^n$. Denote by $a (n) > 0$ the euclidean volume of the euclidean unit
ball in ${\bf R}^n$. The following is known as the {\sl Blaschke-Santal\'o
inequality}, it can be found in [M-P]. (I am grateful to V. Bangert for this
reference).
\bigskip
{\bf Proposition 3.3:} {\it $\rho^n_0 (B_1 \times B_2) \le a (n)^2$ for every
dual pair $(B_1, B_2)$ of convex bodies in ${\bf R}^n$, with equality if and only
if $(B_1, B_2)$ is an euclidean dual pair.}
\bigskip
View ${\bf R}^n$ as an euclidean vector space, equipped with the inner product $<,
>$. For $r > 0$ let $B(r)$ be the closed ball of radius $r$ about $0$ in
${\bf R}^n$; then
$B(r) \times B(r)= D (r)$
is a neighborhood of $0$ in
${\bf R}^{2n}= {\bf R}^n \times {\bf R}^n$.
For $c \ge 1$ define a {\it $c$-quasisymplectic $r$-ball} in ${\bf R}^{2n}$ to be
the image of $D(r)$ under a map $\beta \colon D (r) \to {\bf R}^{2n}$ with the
following properties:

%{\setindent{2)}
{
\item{1)} There are continuous maps $\beta_i\colon B (r) \to {\bf R}^n (i = 1,
2)$ such that $\beta = (\beta_1, \beta_2).$
\item{2)} $\mid \v_0 (x, y) - \v_0 (\beta_1 x, \beta_2 y) \mid \le (c-1) r^2$
for all $x, y \in B(r)$.

}

The above discussion indicates that a $1$-quasisymplectic $1$-ball in
${\bf R}^{2n}$ is an {\sl euclidean dual pair}, i.e. there is a linear isomorphism
$L
\colon {\bf R}^n \to {\bf R}^n$ such that $\beta_1 = L$ and $\beta_2 = (L^t)^{-1}$
where $L^t$ is the transpose of $L$ with respect to the duality ${\bf R}^n \to
({\bf R}^n)^*$ defined by the symplectic form $\rho_0$.

The next lemma contains a more precise statement of the relation between
quasisymplectic balls and euclidean dual pairs.
\bigskip

{\bf Proposition 3.4:} {\it For every $n \ge 1$ there is a number $q (n) > 0$
with the following property:
For every $\ep \in (0, 2^{-n})$, every $r > 0$ and every $(1+
\ep)-$quasisymplectic
$r$-ball
$D$ in ${\bf R}^{2n}$ there is an euclidean dual pair $(B_1, B_2) \subset {\bf R}^n
\times {\bf R}^n$ such that $(1 + \ep^{1/2^n} q (n))^{-1} r (B_1, B_2) \subset D
\subset
(1+ \ep^{1/2 n} q (n)) r (B_1, B_2).$
\bigskip
\it Proof:}
The statement of the proposition is invariant under
scaling, so it is enough to show the proposition for
quasisymplectic 1-balls. For this we proceed by induction on $n = \dim
{\bf R}^{2n}/2.$

The case $n = 1$ is immediate with $q (1) = 1.$ Thus assume that the
proposition is known for $n - 1 \ge 1$ and let $\beta \colon D (1) \to
{\bf R}^{2n}$ be a $(1+\ep)$-quasisymplectic 1-ball where $\ep \le 2^{-n}$. Write
$D_1 = \beta_1 (B (1)), D_2 = \beta_2 (B(1)).$

Let $x_1 \in \partial B (1)$ be a point from the boundary of $B (1)$ and write
$y_1 = \beta_1 (x_1)$. Let $x_2 \in \partial B (1)$ be the unique point which
satisfies $\v_0 (x_1, x_2) = 1$ and write $y_2 = \beta_2 (x_2)$. By our
assumption we have $1 + \ep \ge \v_0 (y_1, y_2) \ge 1 - \ep.$

Define $H_1 = \{ z \in {\bf R}^n \mid \v_0 (z, y_2) = 0\}, H_2 = \{ z \in {\bf R}^n
\mid \v_0 (y_1,z) = 0\}$ and let $P_i \colon {\bf R}^n
={\bf R}\oplus H_i \to H_i$ be the canonical
projection. Define
$W_1 = \{ x \in {\bf R}^n \mid \v_0 (x, x_2) = 0\},
W_2 = \{ x \in {\bf R}^n \mid \v_0 (x_1, x) = 0\}$ and $\bar\beta_i = P_i \circ
\beta_i \mid W_i$. Then
$(\bar\beta_1, \bar\beta_2)$ is a map of the symplectic vector
space $(W_1 \times W_2, \rho_0)$ into the symplectic vector space
$(H_1 \times H_2, \rho_0)$. Moreover the restriction of $\rho_0$ to $W_1
\times W_2$ is the standard symplectic form on
${\bf R}^{2n - 2}$, and the same is true for
$\rho_0\mid_{H_1 \times H_2}.$

We claim that $(\bar\beta_1, \bar\beta_2)$ is a
$1 + \ep +2\ep^2$-quasisymplectic 1-ball in
${\bf R}^{2n - 2}$.

To see this let $z_1 \in W_1$ and write
$\beta_1 (z_1) = a_1 y_1 + \bar\beta_1 (z_1)$. Then $\v_0 (z_1, x_2) = 0$ and

$$\vert\v_0 (\beta_1 z_1, \beta_2 x_2)\vert = \vert\v_0 (a_1 y_1 + \bar\beta_1
(z_1), y_2)\vert \ge\vert a_1\vert (1 - \ep).$$

Since $\beta$ is an $(1 + \ep)$-quasisymplectic 1-ball we conclude that $\mid
a_1 \mid \le \ep/1 - \ep$. The same argument applied to $z_2 \in W_2$ shows
that
$\beta_2 (z_2) = a_2 y_2 + \bar\beta_2 (z_2)$ with $\mid a_2 \mid \le \ep/1
-\ep$. Then
$$\mid \v_0 (\bar\beta_1 (z_1), \bar\beta_2 (z_2)) -  \v_0 (\beta_1
(z_1), \beta_2 (z_2))\mid = \mid a_1 a_2 \v_0 (y_1, y_2)\mid \le \ep^2 (1 +
\ep) /(1 - \ep)^2$$
and hence
$$\mid\v_0 (z_1,z_2) -  \v_0 (\bar\beta_1 (z_1), \bar\beta_2
(z_2))\mid \le \ep + \ep^2(1+\ep)/(1 - \ep)^2 \le 
\ep + 2\ep^2 \le 2\ep.$$

We apply now the induction hypothesis to $(\bar\beta_1, \bar\beta_2).$
This means that there is a linear isomorphism $L_1 \colon W_1 \equiv {\bf R}^{n -
1} \to H_1 \equiv{\bf R}^{n - 1}$ such that
$$\eqalign{(1 + (2 \ep)^{1/2^{n-1}} q (n - 1))^{-1} &(L_1 B(1), (L_1^t)^{-1} B
(1)) \subset (\bar\beta_1 B(1), \bar\beta_2 B(1))\cr
 \subset (1+(2 \ep)^{1/2^{n-1}} q (n - 1) ) &(L_1 B (1), (L_1^t)^{-1} B
(1)).}$$

Define a linear map $L \colon {\bf R}^n \to {\bf R}^n$ by $L \mid_{W_1} = L_1$ and
$L x_1 = y_1.$

Let $x = a_1 x_1 + z_1 \in {\bf R} x_1 \op W_1 \cap B (1)$ and write $\beta_1 (x)
= b_1 y_1 + \bar z_1 \in {\bf R} y_1 \op H_1.$ 
Then we have $\max \{ \v_0 (x, w) \mid w
\in W_2\cap B (1) \} = \pl z_1 \pl \le 1$ and $\mid\v_0 (x, w) - \v_0 (\beta_1
x, \bar\beta_2 w) \mid \le 2 \ep$ for all $w \in W_2 \cap B (1)$.
By induction hypothesis,
$\bar\beta_2 (W_2 \cap B (1))\supset
(1 + 2 \ep^{1/2^{n-1}} q (n - 1))^{-1} (L_1^t)^{-1} B (1)$
and from this we conclude that
$\bar z_1= P_1 \beta_1 (x) \in (\pl z_1 \pl + 2\ep)
(1 + 2 \ep^{1/2^{n-1}} q (n - 1))L_1 B (1)-
(\pl z_1 \pl -2\ep)
(1 + 2 \ep^{1/2^{n-1}} q (n - 1))^{-1} L_1 B (1)$.

Similarly, $a_1 = \v_0 (x, x_2)$ and $ \mid \v_0 (x, x_2) - \v_0 (\beta_1 x,
\beta_2 x_2 ) \mid \le \ep$ implies that
$$\mid b_1 \mid = \v_0 (\beta_1 x, \beta_2 x_2 )
\in [a_1 - \ep, a_1 + \ep].$$

Therefore, if $(\pl\; \pl_L) $ is the norm of the image of the euclidean
scalar product under $L$, then
$$\pl \beta_1 x \pl^2_L = b_1^2 + \pl \bar z_1 \pl^2_L \le (a_1 + \ep)^2 +
(\pl z_1 \pl_L +2\ep)^2 (1 + 2 \ep^{1/2^{n-1}} q (n - 1))^2
\le a_1^2 +
\pl z_1 \pl^2 + \ep^{1/2^{n-1}}  p$$
where $p \ge 1$ only depends on $n$. Together we conclude that $\pl \beta_1
x\pl_L  \le
1 + \ep^{1/2^n} \cdot q (n)$ where again $q (n) \ge 1$ is a
universal constant.

In other words, $\beta_1 B(1) \subset
(1 + \ep^{1/2^n} q (n)) B_1$ where
$B_1$ is the unit ball for $(\pl \; \pl_L)$. The same argument for $\beta_2$
then shows that $(D_1, D_2) \subset
(1 + \ep^{1/2^n} q (n)) (B_1, B_2).$

This is the first half of our claim. The second half follows from the same
argument by taking into account that $\pl \beta_1 x\pl^2_L \ge
(1 + \ep^{1/2^{n-1}} q (n))^{-1}$ for all $x \in \partial B (1).$
{\bf q.e.d.}
\bigskip

From Proposition
3.2 the following is immediate:
\bigskip

{\bf Corollary 3.5:} {\it
For $\ep \le 2^{-n} $ the volume of every $(1 + \ep)$-quasisymplectic
$r$-ball in ${\bf R}^n$ is contained in the interval
$[(1 + \ep^{1/2^n} q (n))^{-2n}r^{2n} a (n)^2,
(1 + \ep^{1/2^n} q (n))^{2n}r^{2n} a (n)^2].$ }
\bigskip
Let now again $\partial \tm$ be the ideal boundary of the universal convering
of a negatively curved manifold $M$ of dimension $n$. Denote by [ ] the
logarithm of the cross ratio of the length cocycle of the metric on $M$.

Using the fact that
$\partial\tm \times \partial \tm - \Delta$
is a smooth symplectic manifold and that for every $\xi \in \partial\tm$ the
subsets $\partial \tm \times \{\xi\}$ and
$\{\xi\} \times \partial \tm $ are Lagrangian we can define for every
fixed $(\xi,\eta) \in \partial\tm \times \partial \tm - \Delta$ a function
$\v$ as in Lemma 3.1. We are going to see that this function is in fact given
by the cross ratio $[ \; ]$. First, for a compact set $K \subset
\partial\tm$ with dense interior and non-empty complement define a function
$\rho_K$ on $(\partial \tm - K )^2$ by $\rho_K (\xi, \eta) = \sup \{[\xi,
\eta,
\zeta, \bar\zeta]\mid \zeta, \bar\zeta \in K\}.$

Then we have:
\bigskip
{\bf Lemma 3.6:} {\it $\rho_K$ is a distance on $\partial \tm - K$.
\bigskip
Proof:} First of all we clearly have $\rho_K (\xi, \xi) = 0$ for all $\xi \in
\partial \tm - K$. Second, since
$[\xi, \eta, \zeta, \bar\zeta]
= - [\eta,\xi, \zeta, \bar \zeta] =
[\eta, \xi, \bar \zeta, \zeta]$
the function $\rho_K$
is symmetric and non-negative.

To show the triangle inequality let $\xi, \eta \in \partial \tm - K$ and let
$\zeta, \bar\zeta \in K$ be such that $ [\xi, \eta, \zeta, \bar \zeta]= \rho_K
(\xi, \eta)$; such points exist by continuity of [ ] and compactness of $K$.
If
$\omega \in \partial \tm - K$ is arbitrary, then
$\rho_K (\xi, \eta) =
[\xi, \omega, \zeta, \bar \zeta]+ [ \omega,\eta, \zeta, \bar \zeta]\le \rho_K
(\xi, \omega) + \rho_K (\omega, \eta).$

We are left with showing that
$\rho_K (\xi, \eta) > 0$ for $\xi \ne \eta$; for this we need the fact that
the interior $U$ of $K$ is non-empty. Choose $\zeta \in U$, let $\gamma$
be a geodesic joining
$\gamma (- \infty) = \zeta$ to $\gamma (\infty) = \xi$, and let $\bar\gamma$
be a geodesic joining
$\bar\gamma (- \infty) = \zeta$ to $\bar\gamma (\infty) = \eta$. We assume
that $\gamma, \bar\gamma$ are parametrized in such a way that for every
$t \in {\bf R}, \gamma' (t)$ and
$\bar \gamma' (t)$ lie on the same strong unstable manifold.
Corollary 2.10 of [H3] then shows that there is
$\bar \zeta \in U \subset K$ in an arbitrarily small neighborhood of $\zeta$
such that $[\xi, \eta, \zeta, \bar\zeta] > 0.$ From this the lemma follows.
{\bf q.e.d.}
\bigskip
Recall now from [H3] that $\partial \tm$ admits a $C^1$-structure if there is
a differentiable structure for $\partial \tm$ for which $\pi$ is a
$C^1$-submersion. Then for every $v \in T^1 \tm$ the restriction to $W^{su}
(v)$ of the natural projection $\pi \colon T^1 \tm \to \partial \tm$ is a
$C^1$-diffeomorphism onto $\partial \tm - \pi (-v)$. In particular we can
equip $\partial \tm$ with a continuous Riemannian metric, and every two such
metrics are globally bilipschitz-equivalent. If $\rho$ is any distance
function on $\partial \tm$ which is locally bilipschitz equivalent to one (and
hence every) continuous Riemannian metric on $\partial \tm$ we say that
$\rho$ is {\sl locally equivalent} to a Riemannian structure. Then we have:
\bigskip
{\bf Proposition 3.7:} {\it If $\partial \tm$ admits a $C^1$-structure then
the distance functions $\rho_K$  on $\partial \tm - K$ are locally equivalent
to a Riemannian structure.
\bigskip
Proof:} Let $\emptyset \ne U \subset \partial \tm$ be open and connected with
closure $K$ and let
$\zeta \in U$ be arbitrarily fixed. If $\eta, \xi \in \partial \tm - K =
\Omega$ and if $
\zeta_1, \zeta_2 \in K$ are such that $\rho_K (\xi, \eta) =
[\xi, \eta, \zeta_1, \zeta_2]$ then $[\xi, \eta, \zeta_1, \zeta_2] =
[\xi, \eta, \zeta_1, \zeta] +[\xi, \eta, \zeta_, \zeta_2]$ shows that
$\rho_K (\xi, \eta) \le 2 \sup \{\vert
[\xi, \eta, \zeta, \bar\zeta] \vert\mid\bar \zeta \in K\} = 2 \bar\rho (\xi,
\eta).$

Thus it suffices to show that the function $\bar\rho$ which is defined by the
above equation is locally equivalent to a Riemannian structure.

For this let $\xi \in \Omega$ again be arbitrary and let $\gamma$ be a
geodesic joining $\gamma (-\infty) = \zeta$ to $\gamma (\infty) = \xi$. Let
$d^{ss}$ (or $d^{su}$) be the distances on the leaves of $W^{ss}$ (or
$W^{su}$ ) induced by the restriction of the Sasaki metric and assume that
$\gamma$ is parametrized in such a way that for $v (\xi) = v = \gamma' (0)$ we
have $\sup\{ d^{ss} (v, w) \mid \pi (-w) \in K\} = 1$. Let moreover $r (\xi) >
0$ be the maximum of all numbers $r > 0$ such that $\{w\mid \pi (-w) \in K\}$
contains the ball of radius $r$ about $v$ in $W^{ss} (v)$. Notice that $\xi
\to v (\xi)$ and $\xi\to r (\xi)$ are continuous.

Let $t \to v (t)$ be a curve of class $C^1$ in $W^{su} (v)$ through $v (0) =
v$. For $w \in W^{ss} (v)$ choose a $d^{ss}$-geodesic $\v_w$ joining $\v_w (0)
= v$
to $\v_w (1) = w$. For $t \ge 0$ let $Y (t) \in T_{\v_w (t)} W^{su}$ be such
that $d \pi (Y (t)) = d \pi (v' (0))$. Then $t \to Y  (t)$ is uniformly
continuous and
$${d \over dt} [ \zeta, \pi (-w), \xi, \pi (v (t)) ] \mid_{t = 0} =
\int^1_0 d \omega (\v_w' (s), Y (s)) \; ds $$
(see [H2]). In other words, if $\chi = \max \{\vert \int^1_0 d \omega (\v'_w
(s), Y(s)) \; ds \vert \mid \pi(-w) \in K \} < \infty$ then $\bar\rho (\xi,
\pi v (t)) = t \chi + o (t)$ and hence $\bar\rho$ is locally Lipschitz at
$\xi$ with respect to the projection of $d^{su} \mid_{W^{su} (v)}$ to
$\partial
\tm$. From continuity and the considerations in [H2] and [H3] the proposition
now immediately follows.
{\bf q.e.d.}
\bigskip
{\it Remark }: Proposition 3.7 shows in particular that
$\partial \tm$ admits a $C^1$-structure only if the distances $\rho_U$
($\emptyset \ne  U\subset \partial \tm$ open) are locally equivalent on the
intersections of their domain of definition. Since the distances $\rho_U$ are
defined just by the length cocycle of the metric this shows that obstructions
to existence of a $C^1$-structure on the geometric boundary $\partial \tm$
of $M$ equipped with a fixed metric $g$ 
can be read off immediately from the
cohomology class of the length cocycle. It seems to be reasonable to believe
that local equivalence
of the distances $\rho_U$ implies the existence of a Lipschitz structure on
$\partial \tm$.

Moreover, if $\partial \tm$ has a $C^1$-structure, then
for every $C^1$-curve $c \colon S^1 \to \partial \tm$ the restriction of
$[\;]$ to quadruples of pairwise distinct points in $c (S^1)$ is a signed
Borel measure, i.e. is $\sigma-$additive and finite. It is not  hard to see
that rectifiable curves in $\partial \tm$ with respect to the $C^1$-structure
are characterized by this property.

We apply now the above considerations to maps of
$({\bf R}^{2n-2}, \rho_0)$
into
$\partial \tm \times \partial\tm - \Delta$
equipped with the cross-ratio $[ \;]$, where as before $\partial \tm$ is the
ideal boundary of the universal covering of a compact negatively curved
manifold $M$ of dimension $n$.

For $c \ge 1$ define a {\sl $c$-quasisymplectic map} of $D
(r)$ into
$\partial \tm \times\partial \tm - \Delta$
 to be a continuous map $\beta$ of $D(r)$ onto a compact subset of
$\partial \tm \times \partial\tm - \Delta$
with the following properties:

%{\setindent{2)}
{
\item{1)} There are continuous maps $\beta_i \colon B(r) \to \partial \tm (i =
1,2)$ such that $\beta = (\beta_1, \beta_2).$
\item{2)} $\mid [\beta_1 (x), \beta_1 (0), \beta_2 (y), \beta_2 (0)] - \v_0
(x, y) \mid \le (c - 1) r^2$ for all $ x \in B (r) , y \in B (r)$.

}

Corollary 2.10 of [H3] then shows that for every $\ep > 0$ there is a number
$r
> 0$ such that for every $\delta \le r$ and every $\xi \ne \eta \in \partial
\tm$ there is a $1 + \ep$-quasisymplectic map $\beta$ of $D(\delta)$
into
$\partial \tm \times \partial\tm - \Delta$
such that $\beta (0, 0) = (\xi, \eta).$

We call the image of $D (r)$ under a $(1 + \ep)$-quasisymplectic map
$\beta$ a $(1 + \ep)$ {\sl quasisymplectic $r$-ball} and call $\beta (0)$
the {\sl center} of the quasisymplectic ball.

For $\ep > 0$ denote by $Q (\ep)$ the family of all $(1 +
\ep)-$quasisymplectic $r$-balls for arbitrary $r > 0$ in
$\partial \tm \times \partial\tm - \Delta$.
Then $Q (\ep) \subset Q (\delta)$ if $\ep \le \delta$. For a
$(1 + \ep)$-quasisymplectic $r$-ball $\beta\colon D (r) \to
\partial \tm \times \partial\tm - \Delta$
define
$$\delta (\beta D (r)) = \sup \{ \vert
[ \xi_1, q_1, \xi_2, q_2] \vert \mid (\xi_1, \xi_2) \in \beta (D (r)), (q_1,
q_2) = \beta (0)\}.$$
Clearly $\delta (\beta D (r) )\in [ (1 - \ep ) r^2, (1 +
\ep ) r^2].$

Fix any distance $d$ on $\partial \tm$ which induces the natural topology and
let
$\diam (B)$ be the
$d \times d$-diameter of a set $B \subset
\partial \tm \times \partial\tm.$ For $\ep > 0$ and a Borel-subset $A$ of
$\partial \tm \times \partial\tm - \Delta$
 define
$$\s_\ep (A) = \inf \{\sum^\infty_{i = 1} \delta (Q_i)^{n - 1} a (n -1)^2 \mid
Q_i
\in Q (\ep), \diam (Q_i) \le \ep, A \subset \cup^{\infty}_{i = 1} Q_i\}$$
and let $\s (A) = \limsup_{\ep \to \infty} \s_\ep (A)$.

Denote by $\lambda = \lambda_g$ the Lebesgue-Liouville measure on
$\partial \tm \times \partial\tm - \Delta$ induced by $g$. Then we have:

\bigskip

{\bf Proposition 3.8:} $\s \le \lambda.$
\bigskip
{\it Proof:} Recall that for every $v \in T^1 \tm$ the restriction to $W^{su}
(v)$ and $-W^{ss} (v)$ of the natural projection $\pi \colon T^1 \tm \to
\partial \tm$ is a homeomorphism onto $\partial \tm - \pi (-v)$ and $\partial
\tm - \pi (v)$. Let $g^{ss}$ be the Riemannian metric on $TW^{ss}$ which is
induced by the Riemannian metric on $M$. Let $\omega$ be the canonical contact
form on $T^1 M$ and $T^1 \tm$ and recall that the bundles $TW^{ss}, TW^{su}$
are Lagrangian with respect to the symplectic form $d\omega$. Thus $d \omega$
determines an isomorphism of $TW^{su}$ onto the dual $(TW^{ss})^*$ of
$TW^{ss}$
and hence the metric $g^{ss}$ induces a metric $g^{su}$ on $TW^{su}$. In fact,
for $p \in \tm$ and for $v \in T^1_p \tm$ the fibre $T_v W^{su}$ of $TW^{su}$
at $v$ can be naturally identified with the $g$-orthogonal complement
$v^\perp$ of  $v$ in $T_p \tm$. Then $v^\perp$ has also an identification with
the tangent space $T^v_v$ at $v$ of the sphere $T^1_p \tm$.
From the explicit form of $d \omega$ (see [H3]) we see that with these
identifications the metric $g^{su}$ on $T_v W^{su}$ is just the natural metric
on $T^v_v$, viewed as the tangent space of the unit sphere $T^1_p \tm$.

Clearly $g^{su}$ is uniformly H\"older continuous, but its restriction to the
leaves of $W^{su}$ is not smooth.

However, for each $v \in T^1 \tm$ we can define smooth coordinates $\v_v
\colon B \to W^{su} (v)$ for $W^{su} (v)$ at $v$ where $B$ is the unit ball in
${\bf R}^{n - 1}$ with the following properties:

%{\setindent{2)}
{
\item{1)} $\v_v (0) = v$ and $d \v_v (0) \colon T_0 {\bf R}^{n - 1} \to (T_v
W^{su}, g^{su})$ is an isometry.
\item{2)} If $\v_v^* g^{su}$ denotes the pull-back of
$g^{su}$ under $\v_v$, then $(x, v) \in B\times T^1 \tm \to \v_v^* g^{su} (x)$
is uniformly H\"older continuous.

}

Let moreover $\exp_v \colon T_v W^{ss} \to W^{ss} (v)$ be the exponential map
of the Riemanian metric $g^{ss}$ and denote by $\pl\pl$ the norm of $TW^i$
induced by $g^i (i = ss, su).$

For $v \in T^1 \tm$ define a map
$$E_v \colon T_v W^{ss} \op T_v W^{su} \to
\partial \tm \times \partial \tm$$
by $E_v (X, Y) = (\pi (- \exp_v X), \pi
(\v_v (Y)))$ (here we identify $T_v W^{su}$ and ${\bf R}^{n -1} $ via $d \v_v
(0))$. Then $E_v$ is a homeomorphism of an open neighborhood $U$ of $0$ in
$T_v W^{ss} \op T_v W^{su}$ onto a neighborhood of $(\pi (-v), \pi (v))$ in
$\partial \tm \times\partial \tm - \Delta$ which is absolutely continuous with
respect to the volume element on $T_v W^{ss} \op T_v W^{su}$ induced by the
symplectic form $d\omega$ and the Lebesgue-Liouville current $\lambda$, with
locally uniformly H\"older continuous Jacobian whose value at
$(0,0) \in T_v W^{ss} \op T_v W^{su}$ equals $1$.

By the definition of the map $E_v$ there is for every fixed $\ep > 0$ a number
$\delta = \delta (\ep) > 0$ such that for every $v \in T^1 \tm$, every $X_i
\in T_v W^i$ with $\pl X_i\pl \le \delta (i = ss, su)$ the following is
satisfied (this is Corollary 2.10 of [H3]):

%{\setindent{ii)}
{
\item{i)} The Jacobian of $E_v$ at $X_{ss} + X _{su}$ is contained in $[1 -
\ep, 1 + \ep].$
\item{ii)} $\vert [ \pi (- \exp_v X_{ss}), \pi (-v), \pi (\v_v X_{su}), \pi
(v)] - d\omega (X_{ss}, X_{su})\vert \le \ep \delta^2.$

}

Let now $\delta \le \delta (\ep)$, let $v \in T^1 \tm$ and let $B^i$ be a ball
of radius $r < \delta$ in $(T_v W^i, g^i).$ By i) above we then have $\lambda
(\pi (- \exp_v B^{ss}), \pi \v_v B^{su} ) \in a (n - 1)^2 r^{2n-2} [1-\ep, 1 +
\ep],$ and ii) shows that the restriction of $E_v$ to $B^{ss} \times B^{su}$
is a $(1 + \ep)$-quasisymplectic $r$-ball. Now if $A \subset
\partial \tm \times\partial \tm - \Delta$ is any compact set,
then for every $\ep > 0$ there are points $v_1, \dots v_k \in T^1 \tm$ and
balls $B_j^i \subset T_{v_j} W^i (i = ss, su)$ of radius $r_j \le  \delta
(\ep)$ such that
$$A \subset \cup^k_{j = 1} E_{v_j} (B_j^{ss} \times B_j^{su}) \quad {\rm and}
\quad \sum^k_{j = 1} \lambda (E_{v_j} (B_j^{ss} \times B_j^{su})) \le \lambda
(A) + \ep.$$

But this means $\lambda (A) \ge (1 - \ep) \sum^k_{j = 1} a (n - 1)^2 r_j^{2n -
2}  - \ep$, and consequently $\s_\ep (A) \le
\lambda (A) + \ep/1-\ep.$

Since $\ep > 0$ was arbitrary and $\s_\delta \ge \s_\ep$ for $\delta\le \ep$
we conclude that $\s_\ep (A) \le \lambda (A)$ for all $\ep > 0$ and therefore
$\s \le \lambda$ as claimed.
{\bf q.e.d.}
\bigskip
We assume now for the moment that the Anosov splitting of $TT^1 M$ is of class
$C^1$. Our goal is to study more precisely the measure $\s$..

For this recall the definition of the {\sl Kanai connection} $\nabla$ on
$T^1M$
which is defined as follows: Let $\j$ be the $(1,1)$-tensor field on $T^1 M$
defined by $\j (Y) = - Y$ for $Y \in TW^{ss}, \j Y =Y$ for $Y \in TW^{su}$ and
$\j ({\bf R} X^0) = 0$. The indefinite metric $h$ on $T^1 M$ defined by $h (Y, Z)
= d\omega (Y, \j Z) + \omega \op \omega (Y, Z)$ is of class $C^1$ and hence
there is a unique affine connection $\nabla$ on $T^1 M$ such that

%{\setindent{ii)}
{
\item{i)} $h$ is parallel with respect to $\nabla$
\item{ii)} The torsion of $\nabla$ is $d \omega \otimes X^0$.

}

The bundles $TW^{su}, TW^{ss}$ are invariant under $\nabla$, and the
restriction of $\nabla$ to every leaf of $W^{su}$ or $W^{ss}$ is flat. In
other words, if we denote the lift of $\nabla$ to $T^1 \tm$ by the same
symbol, then for every $v \in T^1 \tm$ there is a basis of $TW^{su}\mid W^{su}
(v)$ consisting of $\nabla$-parallel vector fields.

The space of geodesics in $\tm$ is
$\partial \tm \times\partial \tm - \Delta = T^1 \tm /{\bf R}$ where ${\bf R}$ acts on
$T^1 \tm$ as the geodesic flow. Equipped with the projection of $d\omega$
this is a smooth symplectic manifold.

According to Darboux's theorem, every point $\xi \in
\partial \tm \times\partial \tm - \Delta$
admits a neighborhood $U$ and local coordinates $(x_1, \dots, x_{n -1},
y_1,\dots , y_{n-1})$ such that in these coordinates the symplectic form
$d\omega$ is the standard symplectic form $\sum d x_i \wedge d
y_i$ on ${\bf R}^{2n -2}.$

We want to use the Kanai-connection $\nabla$ to construct particular such
coordinates which are adapted to the product structure of
$\partial \tm \times\partial \tm - \Delta$.

For this let $v \in T^1 \tm$ be arbitrary and let $\exp_v^{su} \colon T_v
W^{su} \to W^{su} (v)$ and
\hfil\break $\exp_v^{ss} \colon T_v W^{ss} (v) \to W^{ss} (v)$
be
the exponential map of the restriction of $\nabla$ to $W^{su} (v)$ and $W^{ss}
(v)$ respectively. Choose open neighborhoods $A^{su}$ of $v$ in $W^{su} (v),
A^{ss}$ of $v$ in $W^{ss} (v)$ such that for every $w \in A^{su}$ and $z
\in A^{ss}$ the intersection $W^{ss} (w) \cap W^u(z)$ consists of a
unique point $[w, z].$ Then $H = \{ [w, z]\mid w \in A^{su}, z \in A^{ss}\}$
is a local hypersurface in $T^1 \tm$ of class $C^1$ transversal to the
geodesic
flow and hence $(H, d\omega)$ is a symplectic manifold which can naturally be
identified with a neighborhood of $(\pi (v), \pi (-v))$ in
$\partial \tm \times\partial \tm - \Delta$.
For each point $z \in H$ the intersection $W^{ss} (z) \cap H$ is an open
neighborhood of $z$ in $W^{ss} (z).$ Let $X_1, \dots, X_{n-1}$ be a basis of
$T_v W^{su}, Y_1,\dots, Y_{n-1}$ be the basis of $T_v W^{ss}$ which is dual
with respect to $d\omega$ (i.e. such that $d\omega (X_i, Y_j) = \delta_{ij})$
and define local coordinates $\Psi \colon {\bf R}^{n -1} \times {\bf R}^{n-1} \to H$
by $\Psi (x_1, \dots, x_{n - 1}, y_1,\dots, y_{n-1}) = [\exp_v^{su} (\sum x_i
X_i),
\exp^{ss}_v (\sum y_j Y_j)].$
\bigskip

Via these coordinates the symplectic structure on $H$ induces a function $\v$
as in Lemma 3.1.
\bigskip
{\bf Lemma 3.9:} {\it Let $x \in A^{su}, y \in A^{ss}$. Then $[\pi (x), \pi
(v), \pi (-y), \pi (-v)] = \v (\Psi^{-1} [x, y]).$
\bigskip
Proof:} Let $\beta\colon [0, 1]^2 \to H$ be a map of class $C^1$ with $\beta
(0, 0) = v, \beta (s, 0) \in A^{su}, \beta (0, t)\in A^{ss}$ and $\beta (s, t)
= [\beta (s, 0), \beta (0, t)]$. For each fixed $s \in [0, 1]$ the curve $t
\to \beta (s, t)$ is contained in a leaf of $W^{ss}$, and
$d\pi (
{\partial\over \partial s} \beta
(s, t))$ is independent of $t \in [0, 1].$ In
other words, ${\partial\over \partial s} \beta$ is parallel along the curves
$t \to \beta (s, t).$

Lemma 1 of [H2] then shows that

$$\eqalign{
{d\over d s} [  \pi  \beta (s, 0), \pi (v), &  \pi (- \beta (0 ,1 )), \pi (-v)
]  =\cr
{d\over d s} ([ \pi  \beta (s, 0), \pi (v), & \pi (- \beta (0 ,1 )), \pi
(-\beta (s, 1))]  +\cr
[\pi  \beta (s, 0), \pi (v), & \pi (- \beta (s ,1 )), \pi (-v) ]) =\cr
\int^1_0 d\omega (  {\partial\over \partial s} & \beta (s,t),
{\partial\over \partial t} \beta (s,t)) \; dt }$$

since $\pi (-\beta (0, 1)) = \pi (-\beta (s, 1))$ for all $s$. From this the
lemma is immediate.
{\bf q.e.d.}
\bigskip

On the other hand, $X_i, Y_j$ can be extended to $\nabla$-parallel vector
fields along $W^{su} (v)$ which we denote by the same symbol. Define local
coordinates $\Psi_0 \colon {\bf R}^{n - 1} \times {\bf R}^{n - 1} \to H$ by

$$
\eqalign{\Psi_0 & (x_1, \dots, x_{n-1}, y_1, \dots y_{n-1}) =\cr
        \exp^{ss} & (\sum_j y_j Y_j (\exp_v^{su} (\sum x_i X_i))).}$$
\bigskip
{\bf Lemma 3.10:} {\it
$d\omega = \sum^{n-1}_{i = 1} dx_i \wedge dy_i$ in the coordinates $\Psi_0$.
\bigskip
Proof:} Let $c \colon (- \ep, \ep) \to W^{su} (v)$ be the integral curve of a
$\nabla$-parallel vector field along $W^{su} (v)$, and for some
$\nabla$-parallel section $Y$ of $TW^{ss}$ over $W^{su} (v)$ let
$\beta (s, t) = \exp^{ss} t Y (c (s))$. We want to show that
$d\omega ({\partial\over \partial s} \beta,
{\partial\over \partial t} \beta)$ is independent of $s$ and $t$. Since
$d\omega$ is $\nabla$-parallel we conclude that this function is constant
along the curve $s \to \beta (s, 0)$. On the other hand, for every fixed $s
\in (-\ep, \ep)$ we have $\nabla_{{\partial\over \partial t}}
{\partial\over \partial t} \beta = 0$ and consequently (formally)
$$\eqalign{
{\partial\over \partial t}d\omega
( {\partial\over \partial s} \beta, {\partial\over \partial t} \beta) & =
d\omega ( {\nabla\over \partial t} {\partial\over \partial s} \beta,
{\partial\over \partial t} \beta)=\cr
d\omega ( {\nabla\over \partial s} {\partial\over \partial t} \beta,
{\partial\over \partial t} \beta) & =
{1\over 2} {\partial\over \partial s} d\omega
( {\partial\over \partial t} \beta, {\partial\over \partial t} \beta)  = 0}$$
(recall that the torsion of $\nabla$ equals $d\omega \otimes X^0$). From this
the lemma is immediate.
{\bf q.e.d.}
\bigskip
In the sequel we call coordinates
$(x_1, \dots, x_{n-1}, y_1, \dots y_{n-1})$ for
$\partial \tm \times \partial \tm-\Delta$
of class $C^1$ with the property
that the symplectic form $d\omega$ has the standard form $d\omega = \sum dx_i
\wedge dy_i$ {\sl Darboux-coordinates}. The lemma of Darboux says that
Darboux coordinates exist always.
The Darboux coordinates $\Psi_0$ on the symplectic manifold $H$ (which we can
canonically identify with an open subset of
$\partial \tm \times \partial \tm-\Delta$
) are particularly well adapted to the product structure of
$\partial \tm \times \partial \tm-\Delta$.
Namely, for every fixed $y \in {\bf R}^{n-1}$ the set
$\Psi_0 (\{y\}\times {\bf R}^{n-1})$
is contained in a leaf of the strong stable foliation.
However, it is not true that
$\Psi_0 ({\bf R}^{n-1} \times \{y\})$
is contained in a leaf of the unstable foliation. For otherwise the manifold
$\partial \tm \times \partial \tm-\Delta$
equipped with the Kanai connection would be flat which is impossible.

Darboux coordinates with the above additonal properties are by no
means unique. We can choose them compatible
with a symplectic submanifold of
$\partial \tm \times \partial \tm-\Delta$
which is a product.

More precisely, let $v \in T^1 \tm$ and let $A \subset W^{su} (v), B \subset
W^{ss} (v)$ be $k$-dimensional embedded submanifolds containing $v$ in their
interior. Assume that the restriction of $d\omega$ to $L = \{ [x,y]\mid x \in
A, y \in B\}$ is non-degenerate. Let $X_1, \dots, X_{n-1}$ be a basis of $T_v
W^{su}$ such that $X_1, \dots, X_k$ is a basis of $T_vA$ and that $X_{k + 1},
\dots, X_{n - 1}$ is a basis of the annihilator of $T_v B$ in the dual of $T_v
W^{ss}$, where we identify this dual space with $T_v W^{su}$ via $d\omega$.
Let $Y_1, \dots, Y_{n - 1}$ be the basis of $T_v W^{ss}$ which is dual to the
basis $X_1, \dots, X_{n - 1}$; then $Y_1, \dots, Y_k$ is basis of $T_v B$.
Extend $Y_1, \dots, Y_{n - 1}$ by parallel transport to vector fields on $T_v
W^{su}$. For every $x \in A$ and every $j \le k$ the vector $Y_j (x)$ is
tangent to $L$ and hence $Y_1, \dots, Y_k\mid A$ defines a trivialization of
the dual of $TA$ which induces a trivialization of $TA$. Similarly, by the
dual procedure, for each $x \in A$ we obtain a trivialization of $[x, B]$.
Using the exponential map of these
trivializations we obtain as before Darboux coordinates for $L$ which can be
extended to Darboux coordinates for $H$ (at least in a neighborhood of $v$).
\bigskip
{\bf Lemma 3.11:} {\it If $\partial \tm$ has a $C^1$-structure, then for every
$\ep > 0$ there is a number $\chi = \chi (\ep) > 0$ with the following
property: For every $r > 0$ and every $(1 + \chi)$-quasisymplectic map
$\beta$
with diam $\beta (D (r)) < \chi$ there is a symplectic map $\beta_0\colon D (r
(1 + \ep)) \to
\partial \tm \times \partial \tm-\Delta$
such that $\beta_0 D (r (1 -\ep)) \subset
\beta D (r) \subset\beta_0 D (r(1 +\ep)).$
\bigskip
Proof:} The statement of the lemma is local and will be proved by a local
argument.

Let $\rho$ be a continuous symplectic form on ${\bf R}^n \times {\bf R}^n$ with the
following properties:

%{\setindent{5)}
{
\item{1)} At every point in $\{0\} \times {\bf R}^n \cup{\bf R}^n  \times \{0\}, \rho$
coincides with the standard symplectic form $\rho_0 = \sum dx_i \wedge dy_i$.
\item{2)} For every $x \in {\bf R}^n$ the submanifolds $\{x \} \times {\bf R}^n$ and
${\bf R}^n \times \{x\}$ are Lagrangian with respect to $\rho$, in particular we
can define the function $\v$ as in the beginning of this section.
\item{3)} For every point $0 \ne x \in {\bf R}^n$ with sufficiently small
euclidean
norm the sets $A (x) = \{ y \in {\bf R}^n \mid \v (x, y) = 0\}$ and $C (x) = \{y
\in {\bf R}^n \mid \v (y, x) = 0\}$ are $C^1$- hypersurfaces in ${\bf R}^n$ whose
tangent spaces at $0$ coincide with the tangent space at $0$ of the
hypersurface $A_0 (x) = \{ y\mid \v_0 (x, y) = 0\}$ and $C_0 (x) = \{y \in
{\bf R}^n \mid \v_0 (y, x) = 0\}$ where as before, $\v_0$ is defined by the
standard form $\rho_0$.
\item{4)} For a fixed compact ball $B$ about $0$ in ${\bf R}^n$ and a unit vector
$x \in {\bf R}^n$ the hypersurfaces $A (tx) \cap B$ and $C (tx) \cap B$ converge
as $t \to 0$ in the $C^1$-topology to $A_0 (x).$
\item{5)} For every $\ep > 0$ there is a number $\chi = \chi (\ep) > 0$ such
that $\vert \v (x, y) \vert \le (1 + \ep) \parallel x
\parallel \parallel y \parallel$ and $\v (sx, tx) \ge (1 + \ep)^{-1} st \pl
x\pl^2$ for all $x, y \in {\bf R}^n$ with $\pl x\pl < \chi, \pl y\pl < \chi$ and
$s, t \in [0, 1].$

}

Using the explicit form of the local coordinates on
$\partial \tm \times \partial \tm-\Delta$
given in Lemma 3.9 we see that the above properties are satisfied for the
function $\v$ defined by the cross ratio near a given point.

Use the function $\v$ to define a notion of a $(1 + \ep)$-quasisymplectic
$r$-ball on ${\bf R}^n \times {\bf R}^n$. We want to show by induction on $n$ that the
analogue of the statement of the lemma for this function $\v$ is satisfied.
The argument used is similar to the one given in the proof of Proposition
3.4.

The case $n = 1$ is obvious, so assume that the above claim is known for some
$n - 1 \ge 1$.

Let $\v \colon {\bf R}^n \times {\bf R}^n \to {\bf R}$ be as above. By 3) above there is a
number $\chi_1 > 0$ such that for every $0 \ne x \in {\bf R}^n$ with
$\parallel x \parallel < \chi_1$
the euclidean angle between every $0 \ne y \in A (x) \cap B (\chi_1)$ and $A_0
(x)$ is smaller than $\pi/16$.

Let $\kappa \in (0, \pi/16)$, let $\chi = \min\{ \chi_1, \chi(\kappa)\}$ where
$\chi (\kappa) > 0$ is as in 5) above, let
$r >0$ and let $\beta\colon D (r) \to B
(\chi) \times B(\chi)$ be a $(1 + \kappa)$-quasisymplectic embedding. By our
definition we can write $\beta = (\beta_1, \beta_2)$ where $\beta_i$ is a
continuous map of $B (r)$ into ${\bf R}^n$ with $\beta_i (0) = 0$.

For every $x$ from the boundary $\partial B (r)$ of $B (r)$ choose a point $f
(x) \in \partial B (r)$ such that $\beta_2 f (x)$ is contained in the
half-line through $0$ spanned by $\beta_1 (x)$. Such a point exists since
$\beta_2 B (r)$ is a neighborhood of $0$. Since the maps $\beta_i$ are
continuous we can arrange in such a way that $f (x)$ depends continuously
on $x$.

We want to find some $x \in \partial B (r)$ such that the angle between $x$
and $f (x)$ is not larger than $\pi /8$. For this choose a point $x \in
\partial B (r)$ such that the euclidean norm
$\nu_1 = \parallel \beta_1 (x) \parallel$ of $\beta_1 (x)$ is maximal.
For
$\nu_2 = \parallel \beta_2 (fx) \parallel$
we then have $\v (\beta_1 (x), \beta_2 (fx)) / \nu_1\nu_2 \ge (1 +
\kappa)^{-1}$.

By our assumption on a quasisymplectic $r$-ball there is a point $y \in
\partial B (r)$ such that $\v (\beta_1 (y), \beta_2 (fx)) \ge r^2 (1 -
\kappa).$ Since $\v (\beta_1 (y), \beta_2 (fx)) \le (1 + \kappa)
 \parallel \beta_1 (y) \parallel
 \parallel \beta_2 (fx) \parallel$ and
$\pl\beta_1 (y) \parallel \le \nu_1$ we conclude that
$\nu_2 \ge r^2 (1 - \kappa) / (1 + \kappa) \nu_1$ and consequently
$\v (\beta_1 x, \beta_2 (fx)) \ge r^2 (1 - \kappa) / (1 + \kappa)^2$.
Again by the definition of a quasisymplectic $r$-ball this implies that the
angle between $x$ and $fx$ is not larger than $\pi/8$ and is close to zero if
$\chi$ and $\kappa$ are sufficiently close to $0$.

Now let $x$ be a point 
as above; then $C (\beta_2 fx) \times A (\beta_1 x)$ is a
symplectic manifold, and the projection of $\beta \mid C_0 (fx) \times A_0
(x)$ to $C(\beta_2 fx) \times A (\beta_1 x)$ is a $(1 + \kappa')$-
quasisymplectic
embedding where $\kappa' > \kappa$ depends on $\kappa$ and $\chi$ and tends to
zero with $\kappa$ and $\chi$. By our induction hypothesis we can find
symplectic coordinates on $C (\beta_2 fx) \times A (\beta_1 x)$ such that with
respect to these coordinates the projection of $\beta\mid C_0 (f x) \times A_0
(x)$ contains $D (r (1-\ep))$ and is contained in $D (r (1 + \ep))$ where $\ep
> 0$ depends on $\kappa'$
and $\chi$ and tends to zero as $\kappa' \to 0$
and $\chi \to 0$. With the
arguments in the proof of Proposition 3.4 and the explicit way to construct
suitable symplectic coordinates described in Lemma 3.10 we can extend these
symplectic coordinates to coordinates on a neighborhood of $\beta D (r)$ with
the required properties.
{\bf q.e.d.}

\bigskip
As an immediate corollary to Proposition 3.8 and Lemma 3.11 we obtain:

\bigskip
{\bf Corollary 3.12:} {\it If the Anosov splitting of $TT^1 M$ is of class
$C^1$, then $\s = \lambda$.}
\bigskip
We can also use quasisymplectic balls in
$\partial \tm \times \partial \tm - \Delta$ to define a packing
measure as follows: For $\ep >0$ and an open subset $U$ of
$\partial \tm \times \partial \tm - \Delta$ write
$${\cal P}_\ep (U)= \sup \{\sum_{i=1}^\infty \delta(Q_i)^{n-1}
a(n-1)^2\mid Q_i\in Q(\ep), {\rm diam} Q_i \leq \ep,
Q_i\cap Q_j =\emptyset\}$$
and let ${\cal P}(U)= \lim \inf _{\ep \to 0} {\cal P}_\ep (U)$. 
If $A\subset \partial \tm \times \partial \tm -\Delta$ is any
Borel set, then we define ${\cal P}(A)= \inf \{ {\cal P}(U)\mid
U\supset A\}$. From the arguments in the proof of Proposition 3.8 and
Lemma 3.11 we infer:
\bigskip
{\bf Proposition 3.13:} {\it a) $\lambda \leq {\cal P}\quad$
b) If the Anosov splitting of $TT^1M$ is of class $C^1$ then
$\lambda = {\cal P}$.}
\bigskip
{\bf Corollary 3.14:} {\it Let $M$, $N$ be compact negatively curved
manifolds which have the same marked length spectrum. If the Anosov
splitting of $TT^1M$ is of class $C^1$, then $M$ and $N$ have the
same volume.}
\bigskip
{\it Remark:} A more natural way to obtain a current from a
positive H\"older class would be to drop in the definition of
${\cal S}$ the requirement that the covering family consists of
quasisymplectic balls. More precisely, if we define
$$\bar {\cal S}_\ep (A)=
\inf \{\sum_{i=1}^\infty a(n-1)^2
\delta(A_i\times B_i)^{n-1}\mid
A\subset \cup_{i=1}^\infty A_i\times B_i,
\ {\rm diam}\ (A_i\times B_i)\leq \ep\}$$
and $\bar {\cal S}(A)= \lim \sup_{\ep \to 0}\bar {\cal S}_\ep (A)$, then
clearly $\bar {\cal S}\leq {\cal S}$. The Blaschke-Santal\'o
inequality for ${\bf R}^n$ indicates that probably
$\bar {\cal S}={\cal S}$ always. The only case where I was able
to check this 
equality is the case of a hyperbolic 3-manifold. 
 
\bigskip
\bigskip

\centerline{\bf 4. Intersection in higher dimensions}
\bigskip
\bigskip
In this section we consider again a closed negatively curved manifold $M$. As
before, we denote by $\partial \tm$ the ideal boundary of the universal
covering $\tm$ of $M$.

The Riemannian metric $g$ on $M$ lifts to a Riemannian metric $g^u$ on the
leaves of $W^u$. This metric defines a family of Lebesgue measures
$\lambda^{su}$ on the leaves of $W^{su}$. The measures $\lambda^{su}$ are
quasi-invariant under the action of the geodesic flow and they transform via
$${d \over dt} \lambda^{su} \circ \Phi^t \mid_{t = 0} = tr U.$$
Here for every $v \in T^1 M,\; tr U (v)$ is just the trace of the second
fundamental form at
$Pv$ of the horosphere $PW^{su} (v).$

The function $v \to tr U (v)$ is H\"older continuous, positive and of pressure
zero and consequently it defines a positive H\"older class which does not depend
on the particular choice of the family of smooth measures on strong unstable
manifolds. Its Gibbs equilibrium state is just the Lebesgue Liouville measure
$\lambda$ on $T^1 M$. We call the class defined by $tr U$ the {\sl
Lebesgue-class}.

Similarly, the image of $\lambda^{su}$ under the flip $\f$ is a family
$\lambda^{ss}$ of Lebesgue measures on strong stable manifolds. If $U (v)$
denotes the second fundamental operator of the horosphere $PW^{su} (v)$,
viewed as a linear automorphism of the orthogonal complement $v^{\perp}$ of
$v$, then $f (v) = \det (U (v) + U (-v))$ is a positive H\"older continuous
function on $T^1 M$, and $\lambda = dt \times d\lambda^{su} \times f d
\lambda^{ss}$ where $dt$ is the 1-dimensional Lebesgue measure on the
flow-lines of the geodesic flow (see [H3]).

The measure $\lambda$ is invariant under the flip $\f\colon v \to - v$ and
hence the cocycle
$\zeta$ defined by $tr U$ is equivalent to the cocycle $\zeta \circ \f$
defined by $tr U \circ \f$ where an equivalence is given by the function
$\log f$.

Following [H5], we use the cocycle $\zeta \circ \f$ to define a H\"older
continuous kernel $k\colon \hfil\break
\tm \times \tm \times \partial \tm \to {\bf R}$. For
every
fixed $\xi \i \partial \tm$ the function $(x, y) \in \tm \to k (x, y, \xi)$
is smooth. Moreover $k (x,x, \xi) = 0$, and if $v \in T_x^1 \tm$ is such that
$\pi (v) = \xi$, then there is a smooth function $\beta \colon W^{ss} (v) \to
{\bf R}$ such that $k (x, \cdot, \xi )^{-1} (0) =\{ P \Phi^{\beta (w)} w \mid w
\in W^{ss} (v)\}$.
From the kernel $k$ in turn we obtain the Lebesgue cross ratio (see [H5]).
For this fix a point $x \in \tm$ and for $v \ne w \in T_x^1 \tm$ choose a
geodesic $\gamma$ in
$\tm$ joining $\gamma (- \infty) = \pi (w)$ to $\gamma (\infty) = \pi (v)$
and write $\alpha (v, w) = {1 \over 2} (\log f (\gamma' (0)) - k
(\gamma (0), x, \pi (v)) - k (\gamma (0), x, \pi (w)).$
Then $\alpha (v, w)$ does not depend on the choice of $\gamma$. Moreover
by Lemma 1.2 of [H5]
and the above considerations, for fixed $\xi \ne \eta \in \partial \tm$ the
function $x \in \tm \to \alpha (\pi^{-1} (\xi) \cap T_x^1 \tm, \pi^{-1} (\eta)
\cap T_x^1 \tm)$ is smooth.

The {\sl Lebesgue-cross ratio} is then the function [ ] on the space of
quadruples
of pairwise distinct points in $\partial \tm$ defined as follows: Choose $x
\in \tm$, let $v_i \in T_x^1 \tm$ be such that $\pi (v_i) = \xi_i$ and write
$[ \xi_1, \xi_2, \xi_3, \xi_4]=
\alpha (v_1, v_3) +\alpha (v_2, v_4)
- \alpha (v_1, v_4)
- \alpha (v_2, v_3).$

This function admits a continuous extension to the space of quadruples
$( \xi_1, \xi_2, \xi_3, \xi_4)$ with
$\xi_2 \ne \xi_3$ and $\xi_1 \ne \xi_4$ and hence can be viewed as a function
on the space of pairs of oriented geodesics in $\tm$, where we identify
$( \xi_1, \xi_2, \xi_3, \xi_4)$ with the pair
$(\gamma_1, \gamma_2)$ of
geodesics with endpoints
$\gamma_1 (- \infty) =  \xi_3, \gamma_1(\infty)= \xi_1,
\gamma_2 (- \infty) =  \xi_4$ and $\gamma_2(\infty)= \xi_2$.
Since the space
$\g\tm$ of geodesics in $\tm$ is naturally a smooth manifold (via the
identification
$\g\tm = T^1 \tm/{\bf R}$ where ${\bf R}$ acts as the geodesic
flow), the cross ratio is therefore a function on a smooth manifold.
Notice that the Lebesgue cross ratio is in general
different for the cross ration which we used in Section 3.

In the next lemma we investigate the regularity of [ ].

\bigskip
{\bf Lemma 4.1:} {\it If the Anosov splitting of $TT^1M$ is of class $C^k$ for
some $k \ge 1$, then $[\ ]$ is of class $C^k$.}
\bigskip
{\it Proof:} If the Anosov splitting is of class $C^k$, then the functions $tr
U$ and $v \to \det (U(v) + U (-v))$ are of class $C^k$ and $\partial \tm$
admits a $C^k$-structure in the sense of [H3]. Then the kernel
$k \colon \tm \times \tm \times\partial \tm \to {\bf R}$
is of class $C^k$ as well, and from this and the definition of [ ] the lemma
easily follows.
{\bf q.e.d.}
\bigskip
The following proposition is an analogue of Corollary 2.13 of [H3]:
\bigskip

{\bf Proposition 4.2:} {\it If 
$[\ ]$ is of class $C^k$ for some $k \ge 1$, then
$\partial \tm$ admits a $C^k$-structure.}
\bigskip
{\it Proof:} Define a cocycle $\beta (v, t)$ for the geodesic flow on $T^1 M$
and $T^1\tm$ by $\beta (v, \zeta \circ \f (v, t)) = t$. Then $\Psi^t (v) =
\Phi^{\beta (v, t)} v$ defines a H\"older continuous times change for $\Phi^t$
which is smooth along the leaves of the stable foliation. The resulting flow
is Anosov in the following sense:

%{\setindent{2)}
{
\item{1)} For every $v \in T^1 M$ the set $\bar W^{ss} (v) = \{w \in T^1 M
\mid d (\Psi^t v, \Psi^t w) \to 0 (t \to \infty) \}$ is a smoothly embedded
submanifold of $W^s (v)$ which can be realized as the graph of a smooth
function on $W^{ss} (v).$
\item{2)} For every $v \in T^1 M$ the set $\bar W^{su} (v) = \{w \in T^1 M
\mid
d (\Psi^t v, \Psi^t w) \to 0 (t \to - \infty)\}$ is a H\"older submanifold of
$W^u
(v)$ which can be realized as the graph of a H\"older continuous function on
$W^{su} (v).$

}

Let now $A_1, A_2$ be open, relative compact subsets of $\partial \tm$ whose
closures do not intersect. Write $\Omega= \{v \in T^1 \tm \mid \pi (v) \in
A_1, \pi (-v) \in A_2\}$. By 1) and 2) above, for all $v, w \in \Omega$ the
intersection $\bar W^{ss} (v) \cap W^u (w)$ consists of a unique point $q (v,
w).$ Let $\bar \sigma (v, w) \in {\bf R}$ be the unique number $\tau \in {\bf R}$ such
that $q (v, w) \in \Psi^\tau \bar W^{su} (w).$ Since $q (v, \Psi^t w) = q (v,
w)$ and $q (\Psi^t v, w) = \Psi^t q (v, w)$ we have $\bar\sigma (\Psi^t v, w)
=
\bar \sigma (v, w) + t, \bar\sigma (v, \Psi^t w) =\bar\sigma (v, w) - t$ and
$\bar \sigma (v, w) + \bar \sigma (w, v) = [\pi (-v), \pi (-w), \pi (v), \pi
( w) ]$ (compare [H3]).

For a fixed point $v \in \Omega$ and $w_1 \in \bar W^{ss} (v), w_2 \in \bar
W^{su}(v)$ we have $\bar\sigma (w_1, w_2) = 0$ and consequently $[\pi (- w_1),
\pi (-w_2), \pi (w_1), \pi(w_2)] = \bar\sigma (w_2, w_1).$

Assume from now on that the cross ratio is a function of class $C^k$ on $\g
\tm \times \g \tm$ for some $k \ge 1.$

We follow the proof of Corollary 2.13 of [H3] and use the cross ratio [ ]
to construct $C^k$-coordinates for the manifolds $\bar W^{ss}$.

First, let $U \subset T^1 \tm$ be a nontrivial open subset. For $v, w \in T^1
\tm$ write $Cr (v, w) = \bar\sigma (v, w) + \bar\sigma (w, v) = [\pi (-v), \pi
(-w), \pi (v), \pi (w) ]$. We determine inductively points
$$w_1, \dots, w_i \in T^1 \tm\hfil\break (1 \le i \le n - 1)$$
and nontrivial open sets $B_1 \supset \dots
\supset B_i$ of $U$ such that for every $w \in B_i$ and every $j \le i$ the
following is satisfied:

%{\setindent{ii)}
{
\item{i)} The restriction of $f_j = Cr (\cdot, w_j)$ to $B_i \cap \bar W^{ss}
(w)$ does not have critical points.
\item{ii)} $\cap^i_{j = 1} f^{-1}_j (f_j (w)) \cap \bar W^{ss} (w) \cap B_i$
is a $C^k$-embedded connected submanifold of $\bar W^{ss} (w) \cap B_i$ of
codimension $i$.

}

Let $0 \le i \le n - 1$ and assume that $B_i$ and $w_1, \dots, w_i$ are
already determined. Let $w \in B$ and let $c\colon (-\ep, \ep) \to \cap^i_{j =
1} f_j^{-1} (f_j (w)) \cap \bar W^{ss} (w) \cap B_i$ be an embedded curve of
class $C^k$ through $c (0) =w$ with nowhere vanishing tangent. Write
$$\xi = \pi (w), \eta_1 = \pi (-w), \eta_2= \pi (-c (\ep/2));$$
the points $\xi, \eta_1, \eta_2 \in \partial \tm$
are pairwise distinct. We claim that there are
$\beta_1 \ne \beta_2 \in \partial \tm - \{\xi, \eta_1, \eta_2\}$ such that
$[\eta_1, \beta_1, \beta_2, \xi] \ne
[\eta_2, \beta_1, \beta_2, \xi]$. For otherwise we conclude by continuity
that
$[\eta_1, \eta_2, \beta, \xi] = 0$ for all $\beta \in \partial \tm - \{
\eta_1, \eta_2, \xi\}$ and consequently also
$[\eta_1, \eta_2, \beta_1, \beta_2] = [\eta_1, \eta_2, \beta_1, \xi] +
[\eta_1, \eta_2, \xi, \beta_2]= 0$ for all $\beta_1, \beta_2$.

Since [ ] is continuous and $\pi_1 (M)-$invariant and since the action of
$\pi_1 (M)$ on the space of geodesics in $\partial \tm$ is minimal this would
mean that [ ] vanishes identically, a contradiction.

Choose $w_{i + 1} \in T^1 \tm$ in such a way that
$$[\eta_1, \pi (- w_{i + 1}), \pi ( w_{i + 1}), \xi] \ne
[\eta_2, \pi (- w_{i + 1}), \pi ( w_{i + 1}), \xi].$$
Then $Cr (w, w_{i + 1})\ne Cr (c (\ep/2), w_{i + 1}).$
Since the function $f_{i + 1} \colon z \to Cr (z, w_{i + 1})$ is of class
$C^k$ it can be used as a coordinate function for the leaves of $\bar W^{ss}$
on an nontrivial open subset $B_{i + 1}$ of $B_i$ (compare the proof of Lemma
2.12 in [H3]).

From this and the arguments in the proof of Corollary 2.13 of [H3] we conclude
the statement of the proposition.
{\bf q.e.d.}

\bigskip
Let now $M, N$ be closed negatively curved manifolds with isomorphic
fundamental groups $\pi_1 (M) = \pi_1 (N) = \Gamma$. Recall that every
isomorphism of $\pi_1 (M)$ onto $\pi_1 (N)$ induces a $\Gamma$-equivariant
homeomorphism $f \colon \partial \tm \to \partial \tn.$ The
Lebesgue Liouville measure on $T^1 M $ (or $T^1 N$) induces a
$\Gamma$-invariant measure class $\lambda_M$ (or $\lambda_N$) on $\partial
\tm$ ( or $\partial N$), and similarly the Bowen Margulis measure induces an
invariant measure class $\mu_M$ (or $\mu_N$). As an immediate consequence of
Proposition 4.2 we obtain:

\bigskip
{\bf Corollary 4.3:} {\it Assume that the Anosov splittings of $T^1 M$ and
$T^1N$ are of class $C^1$. Then every $\Gamma$-equivariant homeomorphism $f
\colon \partial \tm \to \partial \tn$ which satisfies $f \{\lambda_M, \mu_M\}
\cap \{\lambda_N, \mu_N\} \ne \emptyset$ is a $C^1$-diffeomorphism.}
\bigskip
{\it Proof:} Assume for example that $f \lambda_M = \lambda_N$. Then
$f^4$ maps the Lebesgue Liouville cross ratio on $\partial \tm$ to the
Lebesgue Liouville cross ratio on $\partial \tn$ up to a constant factor.
Since the
$C^1$-structures
on $\partial \tm$ and $\partial \tn$ are determined by coordinate functions
constructed out of the cross ratios (compare the proof of Proposition 4.2) we
conclude that $f$ is a $C^1$-diffeomorphism.
{\bf q.e.d.}
\bigskip
Recall that an {\sl orbit-equivalence} of the geodesic flows on $T^1M$ and
$T^1N$ is a homeomorphism
$\Lambda \colon T^1M \to T^1 N$ which maps every
orbit of the geodesic flow on $T^1 M$ order preserving to an orbit of the
geodesic flow on $T^1 N$. Thus for every $v \in T^1 M$ there is a
homeomorphism $\alpha_v \colon {\bf R} \to {\bf R}$ such that
$\Lambda (\Phi^t v) = \Phi^{\alpha_v(t)}\Lambda (v)$
for all $t \in {\bf R}.$ If $M$ and $N$ have
isomorphic fundamental groups, then their geodesic flows are orbit equivalent,
moreover the orbit equivalence can be chosen to be smooth along the flow lines
of the geodesic flow.
\bigskip
{\bf Corollary 4.4:} {\it Under the hypothesis of Corollary 3.3 there is an
orbit equivalence\hfil\break
$\Lambda \colon T^1M \to T^1 N$
of the geodesic flows which is a $C^1$-diffeomorphism.
\bigskip

Proof:} Assume that the hypothesis of Corollary 3.3 is satisfied; then $f
\lambda_M = \lambda_N$ (after normalization). Let $\lambda_M^{su}$ (or
$\lambda_N^{su}$) be the
family of Lebesgue measures on strong unstable manifolds which is induced from
the lift of the Riemannian metric on $M$ (or $N$). For $v \in T^1\tm$ let
$\Lambda(v)$ be the unique vector in $T^1 \tn$ such that $f (\pi (v)) =\pi
(\Lambda(v)), f (\pi (-v)) =\pi (- \Lambda(v))$ and such that the Jacobian of
$f$ at $\pi (v)$ with respect to the projections of
$\lambda^{su}_M \mid_{W^{su} (v)}$ and
$\lambda^{su}_N \mid_{W^{su} (\Lambda(v))}$ equals 1. Then $v \to \Lambda (v)$
is a $\Gamma$-equivariant map which projects to an orbit equivalence of the
geodesic flows on $T^1 M$ and $T^1 N$.

We claim that $\Lambda$ is in fact a $C^1$-diffeomorphism. To see this, recall
that by our assumption the function $tr U$ on $T^1 M$ and $T^1 N$ is of class
$C^1$ and hence for every $v \in T^1 \tm$ the set
$$\bar W^{ss} (v) = \{w \in W^s (v) \vert d (\pi \lambda^{su}_{M}
\mid_{W^{su} (v)})/ d (\pi \lambda^{su}_{M} \vert_{W^{su} (w)} )(\pi v) = 1\}$$
is a $C^1$-submanifold of $W^s (v).$ The
restriction of $\pi \circ \f$ to $\bar W^{ss} (v)$ is a $C^1$-diffeomorphism
of $\bar W^{ss} (v)$ onto $\partial \tm - \pi (v).$ But $\Lambda \mid_{\bar
W^{ss} (v)} = (\pi \circ \f \mid_{\bar W^{ss} (\Lambda v)} )^{-1} \circ f \circ
(\pi \circ \f \mid_{\bar W^{ss} (v)})$ and hence $\Lambda$ maps every stable
manifold in $T^1 M\; C^1$-diffeomorphically onto a stable manifold in $T^1 N$.
Similarly we also see that $\Lambda$ maps every unstable manifold in $T^1 M\;
C^1$-diffeomorphically onto an unstable manifold in $T^1 N$. In short,
$\Lambda$ is a $C^1$-diffeomorphism.
{\bf q.e.d.}
\bigskip
Recall that the geodesic flows on $T^1 M$ and $T^1 N$ are {\sl homothetic} if
there is a homeomorphism $\Lambda \colon T^1 M \to T^1 N$ and a number $a > 0 $
such that $\Lambda (\Phi^t v) = \Phi^{at} (\Lambda v)$
 for all $v \in T^1 M$ and
all $t \in {\bf R}$. In other words, the metric on $N$ can be rescaled in such a
way that after this rescaling the geodesic flows on $T^1 M$ and $T^1 N$ are
time preserving conjugate.

Next we observe:
\bigskip
{\bf Lemma 4.5:} {\it Let $M, N$ be compact negatively curved manifolds.
Assume that the Anosov splitting of $T^1 M$ and $T^1 N$ is of class $C^1$. If
the geodesic flows on $T^1 M$ and $T^1 N$ are orbit equivalent with an orbit
equivalence of class $C^1$, then they are homothetic.
\bigskip
Proof:} By our assumption, there is a $\Gamma = \pi_1 (M) = \pi_1
(N)$-equivariant diffeomorphism $f \colon \partial \tm \to \partial
\tn$ of class $C^1$,
 and $f \times f$ is a $C^1$-diffeomorphism of the space $\g \tm$ of
geodesics in $\tm$ onto the space $\g \tn$ of geodesics in $\tn$.

The symplectic form $d \omega_0$ on $T^1 \tn$ is invariant under the geodesic
flow and projects to a smooth symplectic form $\eta_0$ on $\g \tn$. Then $(f
\times f)^* \eta_0$ is a continuous, $\pi_1 (M)$-invariant non-degenerate
2-form on $\g \tm$ which can naturally be pulled back to a continuous,
$\pi_1 (M)$ invariant 2-form $\eta$ on $T^1 \tm$ which is invariant under the
geodesic flow.

Now if $\Lambda \colon T^1 M \to T^1 N$ is an orbit equivalence of class $C^1$
induced as above by $f$,
then naturally $\eta = \Lambda^* d\omega_0$ and hence the continuous form
$\eta$ is exact in the sense of distributions. By Theorem A of [H4], this
means
that necessarily $\eta = rd\omega$ for some $r \ne 0$, where $d \omega$ is the
symplectic form on $T^1 M.$ The discussion in Section 3 of this paper then
shows that $f \times f$ preserves the (usual) cross ratio on $\partial \tm$
and $\partial \tn$ up to a constant factor and hence the geodesic flows on
$T^1 M$ and $T^1N$ are homothetic.
{\bf q.e.d.}
\bigskip
As a corollary, we obtain Theorem C from the introduction:
\bigskip
{\bf Corollary 4.6:} {\it Let $M, N$ be closed negatively curved manifolds
with isomorphic fundamental groups and Anosov splitting of class $C^1$. Then
the natural boundary map $\partial \tm \to \partial \tn$ preserves the
Lebesgue measure classes on $\partial \tm$ and $\partial \tn$ if and only if
the geodesics flows on $T^1 M$ and $T^1 N$ are homothetic.}
\bigskip
Let $L$ be the set of length cocycles of smooth metrics of negative
curvature on $M$ with Anosov splitting of class $C^1$. By Corollary 4.6 there
is an injective map $\beta$ of $L$ into the space of geodesic currents $\c M$
of
$M$ which associates to a length cocycle $\alpha$ of a metric $g$ the
Lebesgue Liouville current of $g$. For every $\ell \in L$ and every geodesic
current $\gamma$ on $M$ we then can define the {\sl intersection} between
$\gamma$ and $\beta (\ell)$ by $i (\gamma, \beta (\ell)) = \int \ell d
\gamma$.

Notice that this definition is an extension of the definition in the
2-dimensional case, however $i$ is clearly not symmetric, i.e. in general we
have
$i (\beta (\ell), \beta (\ell')) \ne i (\beta (\ell'), \beta (\ell))$ for
$\ell, \ell' \in L$.

We conjecture that this intersection can be extended to a continuous function
on the product of the space of currents on $M$ with the subspace of Gibbs
currents, in particular Corollary 4.6 should hold without any regularity
assumptions.
\bigskip
\bigskip

\centerline{\bf References}
\bigskip
\bigskip
%{\setindent{[BBBBBB]}
{
\item{[B-C-G]} G. Besson, G. Courtois, S. Gallot, {\it "Entropies et 
rigidit\'es
des espaces localement \hfil\break
 sym\'etriques de courbure strictement n\'egative"},
Geometric and Functional Analysis {\bf 5} (1995), 731-799.
\item{[B1]} F. Bonahon, {\it "The geometry of Teichm\"uller space via
geodesic currents"}, Inv. Math. {\bf 92} (1988), 139-162.
\item{[B2]} F. Bonahon, {\it "Bouts des vari\'et\'es
 hyperboliques de dimension 3
"}, Ann. Math. {\bf 124} (1986), 71-158.
\item{[F]} L. Flaminio, {\it "Local entropy rigidity for hyperbolic
manifolds"}, Communications in Analysis and Geometry {\bf 3} (1995), 555-596.
\item{[H1]} U. Hamenst\"adt, {\it "Time preserving conjugacies of geodesic
flows"}, Erg. Th. \& Dyn. Sys. {\bf 9} (1989), 455-464.
\item{[H2]} U. Hamenst\"adt, {\it " Regularity of time preserving conjugacies
for contact Anosov flows with $C^1$-Anosov splitting"},
Erg. Th. \& Dyn. Sys. {\bf 13} (1993), 65-72.
\item{[H3]} U. Hamenst\"adt, {\it "Regularity at infinity of compact
negatively curved manifolds"}, Erg. Th. \& Dyn. Sys. {\bf 14} (1994), 493-514.
\item{[H4]} U. Hamenst\"adt, {\it "Invariant two forms for geodesic flows"},
Math. Ann. {\bf 301} (1995), 677-698.
\item{[H5]} U. Hamenst\"adt, {\it "Cocycles, Hausdorff measures and cross
ratios"}, to appear in Erg. Th. \& Dyn. Sys.
\item{[HC]} Handbook of convex geometry, Academic Press (1992).
\item{[K]} A. Katok, {\it "Entropy and closed geodesics"},
Erg. Th. \& Dyn. Sys. {\bf 2} (1982), 339-365.
\item{[L]} F. Ledrappier, {\it "Structure au bord des 
vari\'et\'es \`a courbure
n\'egative"}, S\'eminaire de th\'eorie spectrale et g\'eom\'etrie {\bf 13} (1994/95),
Institut Fourier, Grenoble, 97-122.
\item{[M-P]} M. Meyer, A. Pajor, {\it "On the Blaschke-Santal\'o inequality"},
Arch. Math. {\bf 55} (1990), 82-93.
\item{[O]} J.P. Otal, {\it "Le spectre marqu\'e des longeurs des surfaces \`a
courbure n\'egative",} Ann. Math {\bf 131} (1990), 151-162.

}
\bigskip
\bigskip
\bigskip
\centerline{Mathematisches Institut der}
\centerline{Universit\"at Bonn}
\centerline{Beringstr. 1}
\centerline{53115 Bonn}
\centerline{Germany}

\centerline{e-mail: ursula@math.uni-bonn.de}

\bye